\documentclass[twoside,12pt]{article}

\usepackage{stackrel}
\usepackage{amssymb}
\usepackage{amsmath}
\usepackage{bbm}
\usepackage{mathrsfs}
\usepackage{xypic}

\makeatletter
\def\hsmash{\relax 
  \ifmmode\def\next{\mathpalette\mathhsm@sh}\else\let\next\makehsm@sh
  \fi\next}
\def\makehsm@sh#1{\setbox\z@\hbox{#1}\finhsm@sh}
\def\mathhsm@sh#1#2{\setbox\z@\hbox{$\m@th#1{#2}$}\finhsm@sh}
\def\finhsm@sh{\wd\z@\z@ \box\z@}
\makeatother

\makeatletter
\gdef\th@mychange{\normalfont\slshape
   \def\@begintheorem##1##2{\item
        [\hskip\labelsep \theorem@headerfont ##2. ##1  \,--\!--\!--\!--  ]}%
 \def\@opargbegintheorem##1##2##3{%
   \item[\hskip\labelsep \theorem@headerfont ##2. ##1\ {\upshape(}##3{\upshape)}. \,-----  ]}}
\makeatother

\RequirePackage{theorem}
\theoremstyle{mychange}

{\theorembodyfont{\rmfamily}\newtheorem{ttt}{}[section]}
{\theorembodyfont{\rmfamily}\newtheorem{nota}[ttt]{Notation.}}
{\theorembodyfont{\rmfamily}\newtheorem{defi}[ttt]{Definition.}}
{\theorembodyfont{\rmfamily}\newtheorem{remark}[ttt]{Remark.}}
{\theorembodyfont{\rmfamily}\newtheorem{rems}[ttt]{Remarks.}}
{\theorembodyfont{\rmfamily}\newtheorem{ex}[ttt]{Example.}}
{\theorembodyfont{\rmfamily}\newtheorem{cau}[ttt]{Caution.}}

{\theorembodyfont{\rmfamily}\newtheorem{fade}[ttt]{Fact-Definition.}}
{\theorembodyfont{\rmfamily}\newtheorem{prode}[ttt]{Proposition-Definition.}}

{\theorembodyfont{\itshape}\newtheorem{lem}[ttt]{Lemma.}}
{\theorembodyfont{\itshape}\newtheorem{prop}[ttt]{Proposition.}}
{\theorembodyfont{\itshape}\newtheorem{coro}[ttt]{Corollary.}}
{\theorembodyfont{\itshape}\newtheorem{theo}[ttt]{Theorem.}}

{\theorembodyfont{\rmfamily}}
{\theorembodyfont{\rmfamily}}
{\theorembodyfont{\rmfamily}\newtheorem{remo}[ttt]{Remark}}
{\theorembodyfont{\rmfamily}\newtheorem{remso}[ttt]{Remarks}}
{\theorembodyfont{\rmfamily}\newtheorem{exo}[ttt]{Example}}

{\theorembodyfont{\itshape}\newtheorem{obso}[ttt]{Observation}}
{\theorembodyfont{\itshape}}
{\theorembodyfont{\itshape}}
{\theorembodyfont{\itshape}}
{\theorembodyfont{\itshape}\newtheorem{coroo}[ttt]{Corollary}}
{\theorembodyfont{\rmfamily}\newtheorem{cono}[ttt]{Construction}}

{\theorembodyfont{\rmfamily}\newtheorem{proo}[ttt]{Proof of Theorem~\ref{999}.}}

\newcounter{abc}
\newenvironment{abc}{\begin{list}{\rm \alph{abc}) }{\usecounter{abc} \leftmargin=0.0pt \labelsep=0.0pt \listparindent=0.0pt \labelwidth=0.0pt \parsep=\smallskipamount \itemsep=0.0pt \topsep=0.0pt \partopsep=\smallskipamount}}{\end{list}}
\newcounter{iii}
\newenvironment{iii}{\begin{list}{\rm \roman{iii}) }{\usecounter{iii} \leftmargin=0.0pt \labelsep=0.0pt \listparindent=0.0pt \labelwidth=0.0pt \parsep=\smallskipamount \itemsep=0.0pt \topsep=0.0pt \partopsep=\smallskipamount}}{\end{list}}
\newcounter{ABC}
\newenvironment{ABC}{\begin{list}{\rm \Alph{ABC}. }{\usecounter{ABC} \leftmargin=0.0pt \labelsep=0.0pt \listparindent=0.0pt \labelwidth=0.0pt \parsep=\smallskipamount \itemsep=0.0pt \topsep=0.0pt \partopsep=\smallskipamount}}{\end{list}}

\renewcommand{\atop}[2]{\genfrac{}{}{0pt}{}{#1}{#2}}

\newcommand{\bP}{\mathop{\text{\bf P}}\nolimits}
\newcommand{\Spec}{\mathop{\text{\rm Spec}}\nolimits}
\newcommand{\rk}{\mathop{\text{\rm rk}}\nolimits}

\newcommand{\Gal}{\mathop{\text{\rm Gal}}\nolimits}
\newcommand{\Br}{\mathop{\text{\rm Br}}\nolimits}
\newcommand{\inv}{\mathop{\text{\rm inv}}\nolimits}
\newcommand{\Val}{\mathop{\text{\rm Val}}\nolimits}
\newcommand{\ev}{\mathop{\text{\rm ev}}\nolimits}
\newcommand{\Pic}{\mathop{\text{\rm Pic}}\nolimits}
\newcommand{\Div}{\mathop{\text{\rm Div}}\nolimits}
\renewcommand{\div}{\mathop{\text{\rm div}}\nolimits}

\newcommand{\Hom}{\mathop{\text{\rm Hom}}\nolimits}
\newcommand{\res}{\mathop{\text{\rm res}}\nolimits}
\newcommand{\cores}{\mathop{\text{\rm cores}}\nolimits}

\newcommand{\cl}{\mathop{\text{\rm cl}}\nolimits}

\newcommand{\A}{\mathop{\text{\bf A}}}

\newcommand{\St}{\mathop{\text{\rm St}}}

\newcommand{\calO}{\mathscr{O}}
\newcommand{\calQ}{\mathscr{Q}}
\newcommand{\calS}{\mathscr{S}}
\newcommand{\calT}{\mathscr{T}}

\newcommand{\bbA}{{\mathbbm A}}

\newcommand{\bbF}{{\mathbbm F}}

\newcommand{\bbQ}{{\mathbbm Q}}
\newcommand{\bbR}{{\mathbbm R}}
\newcommand{\bbZ}{{\mathbbm Z}}

\newcommand{\frakp}{{\mathfrak{p}}}

\newcommand{\pmodulo}[1]{\nobreak\ifinner\mkern8mu\else\mkern18mu\fi
 (\textup{mod}\,\,#1)}

\newcommand{\br}{ }
\newcommand{\brr}{, }

\def\rightend#1#2{{%
 \leavevmode\nobreak\hskip .5em plus 1fil
 \penalty600 \hskip 0pt plus -1filll
 \vadjust{}\nobreak\hskip 0pt plus 1filll%
 #1\parfillskip=#2\relax \par}}

\def\eop{\ifmmode\rule[-22pt]{0pt}{1pt}\ifinner\tag*{$\square$}\else\eqno{\square}\fi\else\rightend{$\square$}{0pt}\fi}

\author{Andreas-Stephan Elsenhans${}^*$ and J\"org Jahnel${}^\ddagger$}

\date{}

\title{On the order three Brauer classes for cubic~surfaces}

\begin{document}

\maketitle

\begin{abstract}
We~describe a method to compute the Brauer-Manin~obstruction for smooth cubic surfaces
over~$\bbQ$
such that
$\Br(S)/\Br(\bbQ)$
is a
\mbox{$3$-group}.
Our~approach is to associate a Brauer class with every ordered triplet of Galois invariant pairs of Steiner~trihedra. We~show that all order~three Brauer classes may be obtained in this~way. To~show the effect of the obstruction, we give explicit~examples.
\end{abstract}

\footnotetext[1]{Mathematisches Institut, Universit\"at Bayreuth, Univ'stra\ss e 30, D-95440 Bayreuth, Germany,\\
{\tt Stephan.Elsenhans@uni-bayreuth.de}, Website:\! {\tt http://www.staff\!.\!uni-bayreuth.de/$\sim$btm216}}

\footnotetext[3]{\mbox{D\'epartement Mathematik, Universit\"at Siegen, Walter-Flex-Str.~3, D-57068 Siegen, Germany,} \\
{\tt jahnel@mathematik.uni-siegen.de}, Website: {\tt http://www.uni-math.gwdg.de/jahnel}}

\footnotetext[1]{The first author was supported in part by the Deutsche Forschungsgemeinschaft (DFG) through a funded research~project.}

\footnotetext{{\em Key words and phrases.} Cubic surface, Steiner trihedron, Triplet, Twisted cubic curve, Weak approximation, Explicit Brauer-Manin obstruction}

\footnotetext{{\em 2000 Mathematics Subject Classification.} Primary 11D25; Secondary 11D85, 14J26, 14J20}

\section{Introduction}

\begin{ttt}
For~cubic surfaces, weak approximation and even the Hasse principle are not always~fulfilled. The~first example of a cubic surface violating the Hasse principle was constructed by Sir Peter Swinnerton-Dyer~\cite{SD1} in~1962. A~series of examples generalizing that of Swinnerton-Dyer is due to L.\,J.~Mordell~\cite{Mo}. An~example of a different sort was given by J.\,W.\,S.~Cassels and  M.\,J.\,T.~Guy~\cite{CG}.

Around~1970, Yu.\,I.~Manin~\cite{Ma} presented a way to explain these examples in a unified~manner. This~is what today is called the Brauer-Manin~obstruction. Manin's~idea is that a non-trivial Brauer class may be responsible for the failure of weak~approximation. 
\end{ttt}

\begin{ttt}
Let~$S$
be a projective variety
over~$\bbQ$
and
$\pi$
be the structural~morphism. Then,~for the Brauer-Manin~obstruction, only the factor
group~$\Br(S)/\pi^*\!\Br(\bbQ)$
of the Grothendieck-Brauer~group is~relevant.
When~$\smash{\Br(S_{\overline\bbQ}) = 0}$,
that group is isomorphic to the Galois cohomology group
$\smash{H^1(\Gal(\overline{\bbQ}/\bbQ), \Pic(S_{\overline\bbQ}))}$.

For~$S$
a smooth cubic surface, a theorem of Sir Peter~Swinnerton-Dyer~\cite{SD3} states that,
for~$\smash{H^1(\Gal(\overline{\bbQ}/\bbQ), \Pic(S_{\overline\bbQ}))}$,
there are only five~possibilities. It~may be isomorphic
to~$0$,
$\bbZ/2\bbZ$,
$\bbZ/2\bbZ \times \bbZ/2\bbZ$,
$\bbZ/3\bbZ$,
or~$\bbZ/3\bbZ \times \bbZ/3\bbZ$.
\end{ttt}

\begin{ttt}
The~effect of the Brauer-Manin~obstruction has been studied by several~authors. For~instance, the examples of Mordell and Cassels-Guy were explained by the Brauer-Manin~obstruction in~\cite{Ma}. For~diagonal cubic surfaces, in general, the computations were carried out by Colliot-Th\'el\`ene, Kanevsky, and Sansuc in~\cite{CKS}. More~recently, we treated the situations that
$\smash{H^1(\Gal(\overline{\bbQ}/\bbQ), \Pic(S_{\overline\bbQ})) \cong \bbZ/2\bbZ}$
and
$\bbZ/2\bbZ \times \bbZ/2\bbZ$~\cite{EJ2}.
It~seems, however, that for only a few of the cases when
$$H^1(\Gal(\overline{\bbQ}/\bbQ), \Pic(S_{\overline\bbQ})) \cong \bbZ/3\bbZ \, ,$$
computations have been done up to~now. The~goal of the present article is to fill this~gap.
\end{ttt}

\begin{remo}[{\rm On the combinatorics of the 27~lines}{}]
On~a smooth cubic
surface~$S$
over an algebraically closed field, there are exactly 27~lines. The~configuration of these was intensively studied by the geometers of the 19th~century. We~will be able to recall only the most important facts from this fascinating part of classical~geometry. The~reader may find more information in~\cite[Chapter~9]{Do}.

Up~to permutation, the intersection matrix is always the same, independently of the surface~considered. Further,~there are exactly 45~planes cutting out three lines of the~surface. These~are called the {\em tritangent~planes.}

The~combinatorics of the tritangent planes is of interest in~itself. For~us, special sets of three tritangent planes, which are called {\em Steiner trihedra,} will be~important. We~will recall them in section~\ref{drei}.

On~the other hand, there are 72 so-called {\em sixers,} sets of six lines that are mutually~skew. Further,~the sixers come in~pairs. Every~sixer has a partner such that each of its lines meets exacly five of the~other. Together~with its partner, a sixer forms a {\em double-six}. There~are 36 double-sixes on a smooth cubic~surface.

Let~us finally mention one more recent~result. The~automorphism group of the configuration of the 27~lines is isomorphic to the Weyl
group~$W(E_6)$~\cite[Theorem~23.9]{Ma}.
\end{remo}

\begin{ttt}
If~$S$
is defined
over~$\bbQ$
then
$\Gal(\overline\bbQ/\bbQ)$
operates on the 27~lines. The~operation always takes place via a subgroup
of~$W(E_6)$.

The~starting point of our investigations is now a somewhat surprising~observation. It~turns out that
$\smash{H^1(\Gal(\overline{\bbQ}/\bbQ), \Pic(S_{\overline\bbQ}))}$
is of order three or nine only in cases when,
on~$S$,
there are three Galois invariant pairs of Steiner~trihedra that are complementary in the sense that together they contain all the 27~lines. Three~pairs of this sort are classically said to form a~triplet.

This~observation reduces the possibilities for the action
of~$\Gal(\overline\bbQ/\bbQ)$
on the 27~lines. Among~the 350 conjugacy classes of subgroups
of~$W(E_6)$,
exactly 140 stabilize a pair of Steiner~trihedra, while only 54 stabilize a~triplet. Exactly~17 of these conjugacy classes lead to a
$\smash{H^1(\Gal(\overline{\bbQ}/\bbQ), \Pic(S_{\overline\bbQ}))}$
of order three or~nine.
\end{ttt}

\begin{remo}[{\rm Overview over the 350 conjugacy classes}{}]
\label{over}
It~is known~\cite{EJ2} that 
$\smash{H^1(\Gal(\overline{\bbQ}/\bbQ), \Pic(S_{\overline\bbQ}))}$
is of order two or four only when there is a Galois-invariant double-six. There~are the following~cases.

\begin{iii}
\item
Exactly~175 of the 350 classes neither stabilize a double-six nor a~triplet. Then, clearly,
$\smash{H^1(\Gal(\overline{\bbQ}/\bbQ), \Pic(S_{\overline\bbQ})) = 0}$.
\item
Among~the other conjugacy classes, 158 stabilize a double-six and 54 stabilize a~triplet, whereas 37 do~both. For~the latter, one has
$\smash{H^1(\Gal(\overline{\bbQ}/\bbQ), \Pic(S_{\overline\bbQ})) = 0}$,
too.
\item
For~the remaining 17 classes stabilizing a triplet,
$\smash{H^1(\Gal(\overline{\bbQ}/\bbQ), \Pic(S_{\overline\bbQ}))}$
is a non-trivial
\mbox{$3$-group}.
\item
Among~the 121 remaining conjugacy classes that stabilize a double-six, 37 even stabilize a~sixer. Those~may be constructed by blowing up a Galois-invariant set
in~$\bP^2$
and, thus, certainly fulfill weak~approximation.

There~are eight further conjugacy classes such that
$\smash{H^1(\Gal(\overline{\bbQ}/\bbQ), \Pic(S_{\overline\bbQ})) = 0}$,
each leading to an orbit structure of
$[2,5,10,10]$
or~$[2,5,5,5,10]$.
For~the remaining 76 classes,
$\smash{H^1(\Gal(\overline{\bbQ}/\bbQ), \Pic(S_{\overline\bbQ}))}$
is a non-trivial
\mbox{$2$-group}.
\end{iii}

\end{remo}


\begin{ttt}
In~this article, we will compute a non-trivial Brauer class for each of the 16~cases such that
$\smash{H^1(\Gal(\overline{\bbQ}/\bbQ), \Pic(S_{\overline\bbQ})) \cong \bbZ/3\bbZ}$.
We~will start with a ``model case''. This~is the maximal
subgroup~$U_t \subset W(E_6)$
stabilizing a~triplet. We~will show that, for every
subgroup~$G \subset W(E_6)$
leading to
$\smash{H^1(G, \Pic(S_{\overline\bbQ}))}$
of order three, the restriction map from
$\smash{H^1(U_t, \Pic(S_{\overline\bbQ}))}$
is~bijective.
\end{ttt}

\begin{ttt}
We~then present a method to explicitly compute the local evaluation maps induced by an element
$c \in \Br(S)/\pi^*\!\Br(\bbQ)$
of order~three. An~advantage of our approach is that it requires, at most, a quadratic extension of the base~field.

To~evaluate, we test whether the value of a certain rational function is a norm from a cyclic cubic~extension. The~construction of the rational function involves beautiful classical geometry, to wit the determination of a twisted cubic curve
on~$S$.

We~show the effect of the Brauer-Manin obstruction in explicit~examples. It~turns out that, unlike the situation described in~\cite{CKS}, where a Brauer class of order three typically excludes two thirds of the adelic points, various fractions are~possible.
\end{ttt}

\begin{ttt}
Up~to conjugation, there is only one subgroup
$U_{tt} \subset W(E_6)$
such that
$\smash{H^1(U_{tt}, \Pic(S_{\overline\bbQ})) \cong \bbZ/3\bbZ \times \bbZ/3\bbZ}$.
This~group of order three was known to Yu.\,I.\,Manin, in 1969~already.

In~a final section, which is purely theoretic in nature, we show that
$U_{tt}$
actually fixes four~triplets. The~corresponding four restriction maps are injections
$\bbZ/3\bbZ \hookrightarrow \bbZ/3\bbZ \times \bbZ/3\bbZ$.
Together,~each of the eight non-zero elements of the right hand side is met exactly~once.
\end{ttt}

\begin{remo}[{\rm The example of Cassels and Guy}{}]
We~asserted above that all the counterexamples to the Hasse principle, known at the time, were explained by Manin~\cite{Ma}, using Brauer~classes. Concerning~the example of Cassels and Guy~\cite{CG}, this is actually only half the~truth. Let~us recall the facts on this example as far as they are relevant for our~purposes.

The~result of Cassels and Guy is that the cubic
surface~$S$
over~$\bbQ$,
given by
$5x^3 + 12y^3 + 9z^3 + 10w^3 = 0$,
violates the Hasse~principle. To~explain this phenomenon, Manin introduces two Brauer classes
$\smash{\alpha_1 \in \Br(S_{\bbQ(\zeta_3,\sqrt[3]{36},\sqrt[3]{90})})}$
and~$\smash{\alpha_2 \in \Br(S_{\bbQ(\zeta_3,\sqrt[3]{36},\sqrt[3]{30})})}$
that are represented (cf.~Definition~\ref{repr}, below) by
$$\textstyle
f_1 := (x + \frac1{15} \sqrt[3]{90^2} y) / (\sqrt[3]{36}y + 3z)
\quad{\rm and}\quad
f_2 := (z + \frac13 \sqrt[3]{30} w) / (\sqrt[3]{36}y + 3z) \, ,$$
respectively. The~splitting field is the field of definition of all 27 lines, which, in both cases, is a cyclic extension of degree~three. It~turned out only 15 years later~\cite{CKS} that
$\Br(S)/\Br(\bbQ) \cong \bbZ/3\bbZ$
and that
$\alpha_1$
and~$\alpha_2$
are, in fact, the restrictions of Brauer classes
$\underline\alpha_1, \underline\alpha_2 \in \Br(S)$
such
that~$\underline\alpha_2 = -\underline\alpha_1$.

We~will show in Proposition~\ref{restr} that the restriction homomorphisms
$\Br(S_k)/\Br(k) \to \Br(S_{k'})/\Br(k'))$
are automatically injective, as soon as
$k'/k$
is a field extension and both sides are non-trivial
\mbox{$3$-groups}.
Further,~in~\ref{prac}, we give a method that is able to describe
$\underline\alpha_1$
and~$\underline\alpha_2$
explicitly, at least
over~$\bbQ(\zeta_3)$.

Finally,~let us recall that Cassels and Guy work with the rational function
$\smash{f := f_1/f_2 = (x + \frac1{15} \sqrt[3]{90^2} y) / (z + \frac13 \sqrt[3]{30} w)}$.
They~proved that, for every
$x \in S(\bbA_\bbQ)$
for that
$f$~is
defined,
$f(x)$
is a non-principal idel
over~$\bbQ(\sqrt[3]{30},\sqrt[3]{90})$.
This~immediately
implies~$S(\bbQ) = \emptyset$.

Interestingly,~$f$
does not represent a Brauer class, at~all. In~fact,
$\smash{S_{\bbQ(\sqrt[3]{30},\sqrt[3]{90})}}$
is isomorphic to the cubic surface, given by
$5x^3 + 5y^3 + 9z^3 + 9w^3 = 0$.
It~is not hard to see that, in such a case,
$\Br(S_{\bbQ(\sqrt[3]{30},\sqrt[3]{90})}) / \Br(\bbQ(\sqrt[3]{30},\sqrt[3]{90})) = 0$.
\end{remo}

\section{The Brauer-Manin obstruction -- Generalities}
\label{zwei}

\begin{ttt}
For~cubic surfaces, all known counterexamples to the Hasse principle or weak approximation are explained by the following~observation.
\end{ttt}

\begin{defi}
Let~$X$
be a projective variety
over~$\bbQ$
and~$\Br(X)$
its Grothen\-dieck-Brauer~group. Then,~we will~call
$$\ev_\nu \colon \Br(X) \times X(\bbQ_\nu) \longrightarrow \bbQ/\bbZ \, , \quad
(\alpha, \xi) \mapsto \inv\nolimits_\nu (\alpha |_\xi)$$
the {\em local evaluation map\/}.
Here,~$\inv_\nu \colon \Br (\bbQ_\nu) \to \bbQ/\bbZ$
(and~$\inv_\infty \colon \Br (\bbR) \to \frac12\bbZ/\bbZ$)
denote the canonical~isomorphisms.
\end{defi}

\begin{obso}[{\rm Manin}{}]
Let\/~$\pi\colon X \to \Spec \bbQ$
be a projective variety
over\/~$\bbQ$.
Choose an element\/
$\alpha \in \Br (X)$.
Then,~every\/
$\bbQ$-rational
point\/
$x \in X(\bbQ)$
gives rise to an adelic point\/
$(x_\nu)_\nu \in X(\A_\bbQ)$
satisfying the~condition
$$\sum_{\nu\in\Val(\bbQ)} \!\!\!\! \ev_\nu (\alpha, {x_\nu}) = 0 \, .$$
\end{obso}

\begin{rems}
\label{standard}
\begin{iii}
\item
It~is obvious that altering
$\alpha \in \Br(X)$
by some Brauer~class
$\pi^*\!\rho$
for
$\rho \in \Br(\bbQ)$
does not change the obstruction defined
by~$\alpha$.
Consequently,~it is only the factor group
$\Br(X) / \pi^*\!\Br(\bbQ)$
that is relevant for the Brauer-Manin~obstruction.
\item
The~local evaluation map
$\ev_\nu \colon \Br(X) \times X(\bbQ_\nu) \to \bbQ/\bbZ$
is continuous in the second variable.
\item
Further,~for every projective
variety~$X$
over~$\bbQ$
and
every~$\alpha \in \Br (X)$,
there~exists a finite
set~$S \subset \Val(\bbQ)$
such that
$\ev (\alpha, \xi) = 0$
for every
$\nu \not\in S$
and~$\xi \in X(\bbQ_\nu)$.

These~facts imply that the Brauer-Manin obstruction, if present, is an obstruction to the principle of weak~approximation.
\end{iii}
\end{rems}

\begin{lem}
\label{75}
Let~$\pi\colon S \to \Spec \bbQ$
be a non-singular cubic~surface. Then, there is a canonical~isomorphism
$$\delta\colon H^1(\Gal(\overline{\bbQ}/\bbQ), \Pic(S_{\overline\bbQ})) \longrightarrow \Br(S) / \pi^*\!\Br(\bbQ)$$
making the diagram
$$
\definemorphism{gleich}\Solid\notip\notip
\diagram
H^1(\Gal(\overline{\bbQ}/\bbQ), \Pic(S_{\overline\bbQ})) \rto^{\;\;\;\;\;\;\;\;\;\;\delta} \dto_d & \Br(S) / \pi^*\!\Br(\bbQ) \ddto^\res \\
H^2(\Gal(\overline{\bbQ}/\bbQ), \overline\bbQ(S)^*/\overline\bbQ^*) \dgleich & \\
H^2(\Gal(\overline{\bbQ}/\bbQ), \overline\bbQ(S)^*)/\pi^*\!\Br(\bbQ) \rto^{\;\;\;\;\;\;\;\;\;\;\;\;\;\inf} & \Br(\bbQ(S))/\pi^*\!\Br(\bbQ)
\enddiagram
$$
commute.
Here,~$d$
is induced by the short exact~sequence
$$0 \to \overline\bbQ(S)^*/\overline\bbQ^* \to \Div(S_{\overline\bbQ}) \to \Pic(S_{\overline\bbQ}) \to 0$$
and the other morphisms are the canonical~ones.\smallskip

\noindent
{\bf Proof.}
{\em
The~equality at the lower left corner comes from the fact~\cite[section~11.4]{Ta} that
$\smash{H^3 (\Gal(\overline{\bbQ}/\bbQ), \overline\bbQ^*) = 0}$.
The~main assertion is~\cite[Lem\-ma~43.1.1]{Ma}.
}
\eop
\end{lem}

\begin{remark}
\label{fin}
The~group
$\smash{H^1(\Gal(\overline{\bbQ}/\bbQ), \Pic(S_{\overline\bbQ}))}$
is always~finite. Hence, by Remark~\ref{standard}.iii), we know that only finitely many primes are relevant for the Brauer-Manin~obstruction.
\end{remark}

\section{Steiner trihedra and the Brauer group}
\label{drei}

\paragraph{\em Steiner trihedra. Triplets.}

\begin{defi}
\begin{iii}
\item
Three~tritangent planes such that no two of them have one of the 27~lines in common are said to be a~{\em trihedron}.
\item
For~a trihedron
$\{ E_1, E_2, E_3 \}$,
a plane
$E$
is called a {\em conjugate plane\/} if each of the lines
$E_1 \cap E$,
$E_2 \cap E$,
and~$E_3 \cap E$
is contained
in~$S$.
\item
A~trihedron may have either no, exactly one, or exactly three conjugate~planes.
Correspondingly, a trihedron is said to be of the {\em first kind,} {\em second kind,} or {\em third~kind}. Trihedra~of the third kind are also called {\em Steiner~trihedra}.
\end{iii}
\end{defi}

\begin{rems}
\begin{iii}
\item
Steiner~trihedra come in~pairs. Actually,~the three planes conjugate to a Steiner trihedron form another Steiner~trihedron.
\item
All~in all, there are 120 pairs of Steiner trihedra on a non-singular cubic~surface. The~automorphism
group~$W(E_6)$
transitively operates on~them.

The~subgroup
of~$W(E_6)$
stabilizing one pair of Steiner trihedra is isomorphic
to~$[(S_3 \times S_3) \rtimes \bbZ/2\bbZ] \times S_3$
of
order~$432$.
This~group operates on the pairs of Steiner trihedra such that the orbits have lengths
$1$,
$2$,
$27$,
$36$,
and~$54$.
The~orbit of size two is of particular~interest.
\end{iii}
\end{rems}

\begin{fade}
For~every pair of Steiner trihedra, there are exactly two other pairs having no line in common with the pair~given. An~ordered triple of pairs of Steiner trihedra obtained in this way will be called a {\em triplet}.
\end{fade}

\begin{rems}
\begin{iii}
\item
Together,~a triplet contains all the 27~lines.
\item
There~are 240 triplets on a non-singular cubic surface, corresponding to the 40~{\em decompositions\/} of the 27~lines into three pairs of Steiner~trihedra. Clearly,~the operation
of~$W(E_6)$
is transitive on triplets and, thus, on~decompositions.

The~largest subgroup
$U_t$
of~$W(E_6)$,
stabilizing a triplet, is of
order~$216$.
It~is a subgroup of index two
in~$[(S_3 \times S_3) \rtimes \bbZ/2\bbZ] \times S_3$.
As~an abstract group,
$U_t \cong S_3 \times S_3 \times S_3$.
\end{iii}
\end{rems}

\begin{nota}
Let~$l_1, \ldots, l_9$
be the nine lines defined by a Steiner~trihedron. Then,~we will denote the corresponding pair of Steiner~trihedra by a rectangular symbol of the form
$$\left[
\begin{array}{ccc}
l_1 & l_2 & l_3 \\
l_4 & l_5 & l_6 \\
l_7 & l_8 & l_9
\end{array}
\right] \! .$$
The~planes of the trihedra contain the lines noticed in the rows and~columns.
\end{nota}\pagebreak[3]

\begin{ttt}
To~describe Steiner trihedra explicitly, one works best in the blown-up model. In~Schl\"afli's notation~\cite[p.\,116]{Sch} (cf.~Hartshorne's notation in~\cite[Theorem~V.4.9]{Ha}), there are 20~pairs of type~I,
$$\left[
\begin{array}{ccc}
a_i & b_j & c_{ij} \\
b_k & c_{jk} & a_j \\
c_{ik} & a_k & b_i
\end{array}
\right],
$$
10 pairs of type~II,
$$\left[
\begin{array}{ccc}
c_{il} & c_{jm} & c_{kn} \\
c_{jn} & c_{kl} & c_{im} \\
c_{km} & c_{in} & c_{jl}
\end{array}
\right],
$$
and 90 pairs of type~III,
$$\left[
\begin{array}{ccc}
a_i & b_l & c_{il} \\
b_k & a_j & c_{jk} \\
c_{ik} & c_{jl} & c_{mn}
\end{array}
\right].
$$
\end{ttt}

\begin{ttt}
Consequently,~there are two types of decompositions of the 27~lines, 10~consist of two pairs of type~I and one pair of type~II,
$$
\left\{
\left[
\begin{array}{ccc}
a_i & b_j & c_{ij} \\
b_k & c_{jk} & a_j \\
c_{ik} & a_k & b_i
\end{array}
\right],
\left[
\begin{array}{ccc}
a_l & b_m & c_{lm} \\
b_n & c_{mn} & a_m \\
c_{ln} & a_n & b_l
\end{array}
\right],
\left[
\begin{array}{ccc}
c_{il} & c_{jm} & c_{kn} \\
c_{jn} & c_{kl} & c_{im} \\
c_{km} & c_{in} & c_{jl}
\end{array}
\right]
\right\},
$$
30~consist of three pairs of type~III,
$$
\left\{
\left[
\begin{array}{ccc}
a_i & b_l & c_{il} \\
b_k & a_j & c_{jk} \\
c_{ik} & c_{jl} & c_{mn}
\end{array}
\right],
\left[
\begin{array}{ccc}
a_k & b_n & c_{kn} \\
b_m & a_l & c_{lm} \\
c_{km} & c_{ln} & c_{ij}
\end{array}
\right],
\left[
\begin{array}{ccc}
a_m & b_j & c_{jm} \\
b_i & a_n & c_{in} \\
c_{im} & c_{jn} & c_{kl}
\end{array}
\right]
\right\}.
$$
We~will denote these decompositions
by~$\St_{(ijk)(lmn)}$
and
$\St_{(ij)(kl)(mn)}$,~respectively.
\end{ttt}

\begin{remark}
Except~possibly for some of the group-theoretic observations, all the facts on Steiner trihedra given here were well known to the geometers of the 19th~century. They~are due to J.\,Steiner~\cite{St}.
\end{remark}

\paragraph{\em The model case.}

\begin{theo}
\label{999}
Let\/~$\pi\colon S \to \Spec\bbQ$
be a non-singular cubic~surface. Suppose~that\/
$S$
has a\/
$\Gal(\overline\bbQ/\bbQ)$-invariant
triplet and an orbit structure of\/
$[9,9,9]$
on the 27~lines.
Then,~$\Br(S)/\pi^*\!\Br(\bbQ) \cong \bbZ/3\bbZ$.
\end{theo}

\begin{remark}
There~is a different orbit structure of
type~$[9,9,9]$,
in which the orbits correspond to trihedra of the first~kind.
In~this case,
$\Br(S)/\pi^*\!\Br(\bbQ) = 0$.
The~maximal subgroup
$G \subset W(E_6)$
stabilizing the orbits is isomorphic to the dihedral group of
order~$18$.
\end{remark}

\begin{proo}
{\em First step.}
Manin's formula.\smallskip

\noindent
We~have the isomorphism
$\smash{\delta\colon H^1(\Gal(\overline{\bbQ}/\bbQ), \Pic(S_{\overline\bbQ})) \to \Br(S) / \pi^*\!\Br(\bbQ)}$.
In~order to explicitly compute
$\smash{H^1 (\Gal(\overline{\bbQ}/\bbQ), \Pic(S_{\overline\bbQ}))}$
as an abstract abelian group, we will use Manin's formula~\cite[Proposition~31.3]{Ma}. This~means the~following.

$\smash{\Pic(S_{\overline\bbQ})}$~is
generated by the 27~lines. The~group of all permutations of the 27~lines respecting the canonical class and the intersection pairing is isomorphic to the Weyl
group~$W(E_6)$
of
order~$51\,840$.
The~group~$\Gal(\overline{\bbQ}/\bbQ)$
operates on the 27~lines via a finite
quotient~$G := \Gal(\overline{\bbQ}/\bbQ)/H$
that is isomorphic to a subgroup
of~$W(E_6)$.
According~to Shapiro's lemma,
$\smash{H^1(G, \Pic(S_{\overline\bbQ})) \cong H^1(\Gal(\overline{\bbQ}/\bbQ), \Pic(S_{\overline\bbQ}))}$.

Further,~$\smash{H^1 (G, \Pic(S_{\overline\bbQ})) \cong \Hom((N\!D \cap D_0) / N\!D_0, \bbQ/\bbZ)}$.
Here,~$D$
is the free abelian group generated by the 27~lines and
$D_0$
is the subgroup of all principal~divisors.
$N \colon D \to D$
denotes the norm map from
\mbox{$G$-modules}
to abelian~groups.\medskip

\noindent
{\em Second step.}
Divisors.\smallskip

\noindent
The~orbits of size nine define three divisors, which we call
$\Delta_1$,
$\Delta_2$,
and~$\Delta_3$.
By~assumption, these are defined by three pairs of Steiner~trihedra forming a~triplet. Thus,~without restriction,
\begin{eqnarray*}
\Delta_1 & = & (a_1) \;+ (a_2) \;+ (a_3) \;+ (b_1) \;+ (b_2) \;+ (b_3) \;\,+ (c_{12}) + (c_{13}) + (c_{23}) \, , \\
\Delta_2 & = & (a_4) \;+ (a_5) \;+ (a_6) \;+ (b_4) \;+ (b_5) \;+ (b_6) \;\,+ (c_{45}) + (c_{46}) + (c_{56}) \, , \\
\Delta_3 & = & (c_{14}) + (c_{15}) + (c_{16}) + (c_{24}) + (c_{25}) + (c_{26}) + (c_{34}) + (c_{35}) + (c_{36}) \, .
\end{eqnarray*}
Put~$g := \#G/9$.
Then,~the norm of a line is either
$D_1$,
$D_2$,
or~$D_3$
for~$D_i := g\Delta_i$.
Hence,
$$N\!D = \{n_1D_1 + n_2D_2 + n_3D_3 \mid n_1, n_2, n_3 \in \bbZ \} \, .$$
Further,~a direct calculation shows
$LD_1 = LD_2 = LD_3 = 3g$
for every
line~$L$
on~$\smash{\Pic(S_{\overline\bbQ})}$.
Therefore,~a
divisor~$n_1D_1 + n_2D_2 + n_3D_3$
is principal if and only
if~$n_1 + n_2 + n_3 = 0$.

Finally,~$D_0$
is generated by all differences of the divisors given by two tritangent~planes. The~three lines on a tritangent plane are either contained in all three of the divisors
$D_1$,
$D_2$,
and~$D_3$
or only in one of~them. Consequently,
$$N\!D_0 = \langle 3D_1 - 3D_2, 3D_1 - 3D_3, D_1 + D_2 - 2D_3 \rangle \, .$$
The~assertion immediately follows from~this.
\eop
\end{proo}

\begin{rems}
\label{cong}
\begin{iii}
\item
A~generator
of~$\Br(S)/\pi^*\!\Br(\bbQ)$
is given by the~homomorphism
$$\textstyle (N\!D \cap D_0) / N\!D_0 \longrightarrow \frac13\bbZ/\bbZ, \quad n_1D_1 + n_2D_2 + n_3D_3 \mapsto (\frac{n_1 - n_2}3 \bmod \bbZ) \, .$$
\item
Observe~that
$\frac{n_1 - n_2}3 \equiv \frac{n_2 - n_3}3 \equiv \frac{n_3 - n_1}3 \pmod \bbZ$
as
$n_1 + n_2 + n_3 = 0$.
\end{iii}
\end{rems}

\begin{defi}
\label{class}
Let~$\pi\colon S \to \Spec\bbQ$
be a smooth cubic surface and
$\calT = (T_1, T_2, T_3)$~be
a Galois invariant~triplet. This~induces a group~homomorphism
$\Gal(\overline\bbQ/\bbQ) \to U_t$,
given by the operation on the 27~lines.

\begin{iii}
\item
Then,~the image of the generator
$(\frac13 \bmod \bbZ)$
under the natural~homomorphism
$$\textstyle \frac13\bbZ/\bbZ \stackrel{\cong}{\longrightarrow} H^1(U_t, \Pic(S_{\overline\bbQ})) \longrightarrow H^1(\Gal(\overline\bbQ/\bbQ), \Pic(S_{\overline\bbQ})) \stackrel[\delta]{\cong}{\longrightarrow} \Br(S)/\pi^*\!\Br(\bbQ)$$
is called {\em the Brauer class associated with\/} the Galois invariant
triplet~$\calT$.
We~denote it
by~$\cl(\calT)$.
\item
This~defines a map
$$\cl\colon \Theta_S^{\Gal(\overline\bbQ/\bbQ)} \longrightarrow \Br(S)/\pi^*\!\Br(\bbQ)$$
from the set of all Galois invariant triplets, which we will call the {\em class map}.
\end{iii}
\end{defi}

\begin{lem}
Let\/~$\pi\colon S \to \Spec\bbQ$
be a smooth cubic surface and\/
$(T_1, T_2, T_3)$~be
a Galois invariant~triplet. Then,
$$\cl(T_1, T_2, T_3) = \cl(T_2, T_3, T_1) = \cl(T_3, T_1, T_2) \, .$$
On~the other hand,
$\cl(T_1, T_2, T_3) = -\cl(T_2, T_1, T_3)$.\smallskip

\noindent
{\bf Proof.}
{\em
This~is a direct consequence of the congruences observed in~\ref{cong}.ii).
}
\eop
\end{lem}

\paragraph{\em The general case of a Galois group stabilizing a~triplet.}

Using~Manin's formula, we computed
$\smash{H^1 (G, \Pic(S_{\overline\bbQ}))}$
for each of the 350~conjugacy classes of subgroups
of~$W(E_6)$.
The~computations in~{\tt gap} took only a few seconds of CPU~time. We~systematized the results in Remark~\ref{over}. Note~that, in particular, we recovered Sir~Peter~Swinnerton-Dyer's result~\cite{SD3} and made the following~observation.




\begin{prop}
\label{exp}
Let\/~$S$
be a non-singular cubic surface
over\/~$\bbQ$.
If\/
$$H^1 (\Gal(\overline\bbQ/\bbQ), \Pic(S_{\overline\bbQ})) = \bbZ/3\bbZ {\rm ~~or~~} \bbZ/3\bbZ \times \bbZ/3\bbZ$$
then,~on\/~$S$,
there is a Galois invariant~triplet.\smallskip

\noindent
{\bf Proof.}
{\em
This~is seen by a case-by-case study using~{\tt gap}.
}
\eop
\end{prop}

\begin{rems}
\label{3oder9}
\begin{iii}
\item
On~the other hand, if there is a Galois invariant triplet
on~$S$
then
$\smash{H^1 (\Gal(\overline\bbQ/\bbQ), \Pic(S_{\overline\bbQ}))}$
is either
$0$,
or~$\bbZ/3\bbZ$
or~$\bbZ/3\bbZ \times \bbZ/3\bbZ$.
\item
Proposition~\ref{exp} immediately provokes the question whether the cohomology classes are always ``the same'' as in the
$[9, 9, 9]$-cases.
I.e.,~of the type
$\cl(\calT)$
for a certain Galois invariant
triplet~$\calT$.
Somewhat~surprisingly, this is indeed the~case.
\end{iii}
\end{rems}

\begin{lem}
\label{criteria}
Let\/~$\calS$
be a non-singular cubic surface over an algebraically closed field,
$H$~a
group of automorphisms of the configuration of the 27~lines, and\/
$H^\prime \subseteq H$
any~subgroup. Each~of the criteria below is sufficient~for
$$\res\colon H^1 (H, \Pic(\calS)) \to H^1 (H^\prime, \Pic(\calS))$$
being an~injection.

\begin{iii}
\item
$H^\prime$~is
of index prime
to\/~$3$
in\/~$H$
and\/
$H^1 (H, \Pic(\calS))$~is
a\/~\mbox{$3$-group.}
\item
$H^\prime$~is
a normal subgroup
in\/~$H$
and\/~$\rk \Pic(\calS)^H = \rk \Pic(\calS)^{H^\prime}$.
\end{iii}\smallskip

\noindent
{\bf Proof.}
{\em
i)
Here,~$\cores \circ \res \colon H^1 (H, \Pic(\calS)) \to H^1 (H, \Pic(\calS))$
is the multiplication by a number prime
to~$3$,
hence a~bijection.\smallskip

\noindent
ii)
The~assumption ensures that
$H/H^\prime$~operates
trivially
on~$\Pic(\calS)^{H^\prime}$.
Hence,
$H^1 \big( H/H^\prime, \Pic(\calS)^{H^\prime} \big) = 0$.
The~inflation-restriction~sequence
$$0 \to H^1 \big( H/H^\prime, \Pic(\calS)^{H^\prime} \big) \to H^1(H, \Pic(\calS)) \to H^1(H^\prime, \Pic(\calS))$$
yields the~assertion.
}
\eop
\end{lem}

\begin{prop}
\label{restr}
Let\/~$\calS$
be a non-singular cubic surface over an algebraically closed field and\/
$U_t$~be
the group of automorphisms of the configuration of the 27~lines that stabilize a~triplet.\smallskip

\noindent
Further,~let\/~$H \subseteq U_t$
be such
that\/~$H^1 (H, \Pic(\calS)) \cong \bbZ/3\bbZ$.
Then,~the restriction
$$\res\colon H^1 (U_t, \Pic(\calS)) \longrightarrow H^1 (H, \Pic(\calS))$$
is~bijective.\smallskip

\noindent
{\bf Proof.}
{\em
This~proof has a computer~part. Using~{\tt gap}, we made the observation that
$\rk \Pic(\calS)^H = 1$
for every
$H \subset W(E_6)$
such that
$H^1 (H, \Pic(\calS)) \cong \bbZ/3\bbZ$.

To~verify the assertion, we will show~injectivity. According~to Lemma~\ref{criteria}.i), it suffices to test this on the
\mbox{$3$-Sylow~subgroups}
of~$H$
and~$U_t$.
But~the
\mbox{$3$-Sylow~subgroup}
$\smash{U_t^{(3)} \cong \bbZ/3\bbZ \times \bbZ/3\bbZ \times \bbZ/3\bbZ}$
is~abelian. Thus,~Lemma~\ref{criteria}.ii) immediately implies the~assertion.
}
\eop
\end{prop}

\begin{coro}
\label{inj}
Let\/~$H^\prime \subseteq H \subseteq U_t$
be~arbitrary. Then,~for the restriction map\/
$\res\colon H^1 (H, \Pic(\calS)) \longrightarrow H^1 (H^\prime, \Pic(\calS))$,
there are the following~limitations.

\begin{iii}
\item
If\/~$H^1 (H, \Pic(\calS)) = 0$
then\/~$H^1 (H^\prime, \Pic(\calS)) = 0$.
\item
If\/~$H^1 (H, \Pic(\calS)) \cong \bbZ/3\bbZ$
and\/~$H^1 (H^\prime, \Pic(\calS)) \neq 0$
then\/
$\res$
is an~injection.
\item
If\/~$H^1 (H, \Pic(\calS)) \cong \bbZ/3\bbZ \times \bbZ/3\bbZ$
then\/
$H^1 (H^\prime, \Pic(\calS)) \cong \bbZ/3\bbZ \times \bbZ/3\bbZ$
or\/~$0$.
In~the former case,
$H^\prime = H$.
In~the latter case,
$H^\prime = 0$.
\end{iii}\smallskip

\noindent
{\bf Proof.}
{\em
We~know from Remark~\ref{3oder9}.i) that both~groups may be only
$0$,
$\bbZ/3\bbZ$,
or
$\bbZ/3\bbZ \times \bbZ/3\bbZ$.\smallskip

\noindent
i)
If~$H^1 (H^\prime, \Pic(\calS))$
were isomorphic
to~$\bbZ/3\bbZ$
or
$\bbZ/3\bbZ \times \bbZ/3\bbZ$
then the restriction from
$U_t$
to~$H^\prime$
would be the zero~map.\smallskip

\noindent
ii) is immediate from the computations~above.\smallskip

\noindent
iii)
This~assertion is obvious as
$H$~is
of order~three.  
}
\eop
\end{coro}

\section{Computing the Brauer-Manin obstruction}

\paragraph{\em A splitting field for the Brauer class.}

\begin{ttt}
The~theory developed above shows that the non-trivial Brauer classes are, in a certain sense, always the same, as long as they are of order~three. This~may certainly be used for explicit~computations. The~Brauer class remains unchanged under suitable restriction~maps. These~correspond to extensions of the base~field.

We~will, however, present a different method~here. Its~advantage is that it avoids large base~fields.
\end{ttt}

\begin{lem}
Let\/~$S$
be a non-singular cubic~surface and let\/
$T_1$
and\/~$T_2$
be two pairs of Steiner trihedra that define disjoint sets of~lines.\smallskip

\noindent
Then,~the 18~lines defined by\/
$T_1$
and\/~$T_2$
contain exactly three double-sixes. These~form a triple of azygetic double-sixes.\medskip

\noindent
{\bf Proof.}
{\em
We~work in the blown-up model.
As~$W(E_6)$
operates transitively on pairs of Steiner trihedra, we may suppose without restriction~that
$$T_1 =
\left[
\begin{array}{ccc}
a_1 & b_2 & c_{12} \\
b_3 & c_{23} & a_2 \\
c_{13} & a_3 & b_1
\end{array}
\right] \, .$$
Further,~the stabilizer of one pair of Steiner trihedra acts transitively on the two complementary~ones. Thus,~still without~restriction,
$$T_2 =
\left[
\begin{array}{ccc}
a_4 & b_5 & c_{45} \\
b_6 & c_{56} & a_5 \\
c_{46} & a_6 & b_4
\end{array}
\right] \, .$$

In~this situation, we immediately see the three double-sixes
$$
\left(
\begin{array}{cccccc}
a_1 & a_2 & a_3 & a_4 & a_5 & a_6 \\
b_1 & b_2 & b_3 & b_4 & b_5 & b_6
\end{array}
\right),
\left(
\begin{array}{cccccc}
a_1 & a_2 & a_3 & c_{56} & c_{46} & c_{45} \\
c_{23} & c_{13} & c_{12} & b_4 & b_5 & b_6
\end{array}
\right),
\left(
\begin{array}{cccccc}
c_{23} & c_{13} & c_{12} & a_4 & a_5 & a_6 \\
b_1 & b_2 & b_3 & c_{56} & c_{46} & c_{45}
\end{array}
\right)
$$
that are contained in the 18~lines.

Finally,~two double-sixes contained within the 18~lines must have at least six lines in~common. It~is well known that, in this case, they have exactly six lines in~common. Two~such double-sixes are classically called azygetic~\cite[Lemma~9.1.4 and the remarks before]{Do}. A~double-six that is azygetic to the first of the three is necessarily of the~form
$$
\left(
\begin{array}{cccccc}
a_i & a_j & a_k & c_{mn} & c_{ln} & c_{lm} \\
c_{jk} & c_{ik} & c_{ij} & b_l & b_m & b_n
\end{array}
\right) .
$$
Obviously,~such a double-six is contained within the 18~lines considered only for
$\{i,j,k\} = \{1,2,3\}$
or~$\{i,j,k\} = \{4,5,6\}$.
}
\eop
\end{lem}

\begin{prode}
Let~$S$
be a non-singular cubic surface over an arbitrary
field~$K$.
Suppose~that
$T_1$
and~$T_2$
are two pairs of Steiner trihedra that are
$\Gal(\overline{K}/K)$-invariant
and define disjoint sets of~lines.

\begin{abc}
\item
Then,~there is a unique minimal Galois extension
$L/K$
such that the three double-sixes, contained within the 18~lines defined by\/
$T_1$
and~$T_2$,
are
$\Gal(\overline{K}/L)$-invariant.
We~will call
$L$
the {\em field of definition\/} of the double-sixes, given by
$T_1$
and~$T_2$.
\item
Actually,~each of the six sixers is
$\Gal(\overline{K}/L)$-invariant.
\item
$\Gal(L/K)$
is a subgroup
of~$S_3$.
\end{abc}\smallskip

\noindent
{\bf Proof.}
a) and~c)
$\Gal(\overline{K}/K)$~permutes
the three double-sixes.
$L$~corresponds
to the kernel of this permutation~representation.\smallskip

\noindent
b)
The~triple intersections of
$T_1$
or~$T_2$
with two of the double-sixes are 
$\Gal(\overline{K}/L)$-invariant.
From~these trios, the sixers may be~combined.
\eop
\end{prode}

\begin{coroo}[{\rm A splitting field}{}]
\label{split}
Let\/~$S$
be a non-singular cubic surface over a
field\/~$K$
with a Galois-invariant
triplet\/~$(T_1,T_2,T_3)$.\smallskip

\noindent
Then,~the field\/
$L$
of definition of the three double-sixes, given by\/
$T_1$
and\/~$T_2$,
is a splitting field for every class
in\/~$\Br(S)/\pi^*\!\Br(K)$.\medskip

\noindent
{\bf Proof.}
{\em
As~$S_L$
has a Galois-invariant sixer, we certainly~\cite[Theorem~42.8]{Ma} have that
$\Br(S_L)/\pi_L^*\!\Br(L) = 0$.
}
\eop
\end{coroo}

\begin{rems}
Let~$S$
be a non-singular cubic surface
over~$\bbQ$
having a Galois-invariant
triplet~$(T_1,T_2,T_3)$.

\begin{iii}
\item
Then,~there are the three splitting fields defined by
$\{T_1, T_2\}$,
$\{T_2, T_3\}$,
and~$\{T_3, T_1\}$.
It~may happen that two of them or all three~coincide. In~general, they are cyclic cubic extensions of a quadratic number~field.
\item
Among the 17 conjugacy classes of subgroups
of~$W(E_6)$
leading to a Brauer class of order three, there are exactly nine allowing a splitting field that is cyclic of degree~three.
\end{iii}
\end{rems}

\begin{remark}
Suppose~that
$S$
is a non-singular cubic surface
over~$\bbQ$
having a Galois-invariant~triplet. Then,~somewhat surprisingly, every splitting field contains one of the three, described as in Corollary~\ref{split}. The~reason for this is the following~observation.
\end{remark}

\begin{lem}
\label{doppelsechs}
Let\/~$G \subset U_t$
be an arbitrary~subgroup. Then,~either\/
$G$
stabilizes a double-six or\/
$\smash{H^1(G,\Pic(S_{\overline\bbQ})) \neq 0}$.\smallskip

\noindent
{\bf Proof.}
{\em
This~follows from the inspection of all subgroups summarized in~\ref{over}.
}
\eop
\end{lem}\pagebreak[3]

\paragraph{\em A rational function representing the Brauer class.}

\begin{ttt}
To~compute the Brauer-Manin obstruction, we will follow the ``informal'' algorithm of Yu.\,I.~Manin~\cite[Sec.~45.2]{Ma}.
For~$L$
a splitting field, we first have, according to Shapiro's lemma,
$\smash{H^1(\Gal(\overline\bbQ/\bbQ), \Pic(S_{\overline\bbQ})) = H^1(\Gal(L/\bbQ), \Pic(S_{\overline\bbQ}))}$.
If~$\Gal(L/\bbQ) \cong S_3$
then a class of order three is preserved under restriction
to~$A_3$.
As~this is a cyclic~group, 
$$H^1(\Gal(L/\bbQ(\sqrt{\Delta})), \Pic(S_{\overline\bbQ})) \cong (\Div_0(S_{\bbQ(\sqrt{\Delta})}) \cap N\! \Div(S_L)) / N\! \Div_0(S_L) \, .$$
\end{ttt}

\begin{defi}
\label{repr}
Let~$S$
be a non-singular cubic surface
over~$\bbQ$
having a Galois-invariant~triplet.
Suppose~that~$\Br(S)/\pi^*\!\Br(\bbQ) \neq 0$.
Fix~a splitting field 
$L \supset \bbQ(\sqrt{\Delta})$
of
degree~$6$
or~$3$
over~$\bbQ$.
Let,~finally,
$\Psi \in \bbQ(\sqrt{\Delta})(S)$
be a rational function such that
$$\div \Psi = N_{L/\bbQ(\sqrt{\Delta})} D$$
for a divisor
$D \in \Div(S_L)$,
but not for a principal~divisor. We~then say that the rational
function~$\Psi$
{\em represents\/} a non-trivial Brauer~class.
\end{defi}

\begin{ttt}
Since~$\Br(S)/\pi^*\!\Br(\bbQ)$
is a
\mbox{$3$-group},
the 27 lines
on~$S_{\bbQ(\sqrt{\Delta})}$
may have the orbit structure
$[3,3,3,3,3,3,3,3,3]$,
$[3,3,3,9,9]$,
or~$[9,9,9]$.
We~distinguish two~cases.\smallskip

\begin{ABC}
\item
There~is an orbit of size~three.

Then,~each such orbit is actually formed by a~triangle.
If~$L$
is chosen such that some of these orbits split then one may put
$\Psi := l_1/l_2$,
where
$l_1$
and~$l_2$
are two linear forms defining distinct split~orbits.

Forms~of the type described have been used before to describe Brauer classes. See~\cite[Sec.\,45]{Ma}, \cite{CKS}, or~\cite{J}.\smallskip
\item
The~orbit structure
is~$[9,9,9]$.

Then,~the Galois-invariant pairs
$T_1, T_2, T_3$
of Steiner~trihedra immediately define the orbits of the~lines. 
We~will work with the splitting
field~$L$
defined
by~$\{T_1,T_2\}$.
Then,~on~$S_L$,
there are three Galois-invariant
sixers~$s_1, s_2, s_3$.
The~Galois group
$G := \Gal(L/\bbQ(\sqrt{\Delta})) \cong \bbZ/3\bbZ$
permutes them~cyclically.
\end{ABC}
\end{ttt}

\begin{cau}
The~three divisors given by the orbits of size nine generate
in~$\Pic(S_{\bbQ(\sqrt{\Delta})})$
a subgroup of index~three. It~turns out impossible to write down a
function~$\Psi$
representing a non-trivial Brauer class such that
$\div \Psi$
is a linear combination of these three~divisors. Thus,~curves~of degree
$> \!1$
need to be~considered.
\end{cau}

\begin{lem}
On every smooth cubic surface, there are 72 two-dimensional families of twisted cubic~curves. Each~such family corresponds to a~sixer. Having~blown down the sixer, the curves correspond to the lines
in\/~$\bP^2$
through none of the six blow-up~points.\smallskip

\noindent
{\bf Proof.}
{\em
This~is a classical result, due to A.~Clebsch~\cite[p.~371/72]{Cl}.
}
\eop
\end{lem}

\begin{cono}[{\rm of function~$\Psi$, theoretical part}{}]
\begin{iii}
\item
Let~$C = C_1$
be a twisted cubic curve
on~$S_L$
corresponding
to~$s_1$.
Then,~the
\mbox{$G$-orbit}
$\{C_1, C_2, C_3 \}$
of~$C_1$
consists of three twisted cubic curves corresponding to the three~sixers.
\item
An~elementary calculation shows that
$C_1 + C_2 + C_3 \sim 3H$
for
$H$
the hyperplane~section.
\item
Further,~$D_1 \sim D_2 \sim D_3 \sim 3H$
for~$D_i$
the sum of the nine lines
corresponding to~$T_i$.
Hence,~$NC - D_1 = C_1 + C_2 + C_3 - D_1$
is a principal~divisor. It~is clearly the norm of a divisor.
\item
We~choose
$\Psi$
such that
$\div(\Psi) = C_1 + C_2 + C_3 - D_1$.
\end{iii}
\end{cono}

\begin{remark}
It~is a tedious calculation to show that
$NC - D_1$
is not the norm of a principal~divisor.
However,~$D_1 - D_2$
is. In~particular,
$NC - D_1$
and
$NC - D_2$
define the same Brauer~class.
\end{remark}

\begin{cono}[{\rm of function~$\Psi$, practical part}{}]
\label{prac}
For~a concrete cubic surface, one may determine a
function~$\Psi$
using the following~strategy.

\begin{iii}
\item
First,~construct a blow-down map
$S_L \to \bP^2_L$
in two steps by first blowing down an orbit of~three. There~is some beautiful classical algebraic geometry involved in this step~\cite[Sec.~9.3.1]{Do}. The~result is a surface
in~$\bP^1 \times \bP^1 \times \bP^1$,
given by a trilinear~form. The~three projective lines appear as the pencils of the planes through three skew lines
on~$S$.
To~give the morphism down
to~$\bP^2$
is then equivalent to bringing the trilinear form into standard~shape.
\item Then,~choose a line through none of the blow-up~points and calculate seven
\mbox{$L$-rational}
points on the corresponding twisted cubic~curve. Generically,~these determine a two-dimensional space of cubic forms
over~$\bbQ(\sqrt{\Delta})$
containing the form that defines the cubic surface
$S_{\bbQ(\sqrt{\Delta})}$~itself.
It~turned out practical to work with a line that connects the preimages of two
\mbox{$\bbQ$-rational}~points.
\end{iii}
\end{cono}\pagebreak[3]

\paragraph{\em Explicit Brauer-Manin~obstruction.}

\begin{lem}
Let\/~$\pi\colon S \to \Spec \bbQ$
be a non-singular cubic surface with a Galois invariant~triplet.
Further,~let\/
$c \in \Br(S)/\pi^*\!\Br(\bbQ)$
be a non-zero~class.\smallskip

\noindent
Let,~finally,
$L$
be a splitting field that is cyclic of degree three over a
field\/~$\bbQ(\sqrt{\Delta})$
and\/
$\Psi$
a rational function
representing\/~$c$.
We~consider the cyclic algebra
$$\textstyle Q := L(S)\{Y\} / (Y^3 - \Psi)$$
over the function
field\/~$\bbQ(\sqrt{\Delta})(S)$.
Here,~$Y t = \sigma(t) Y$
for\/
$t \in L(S)$
and a fixed generator\/
$\sigma \in \Gal(L/\bbQ(\sqrt{\Delta}))$.\smallskip

\noindent
Then,~$Q$
is an Azumaya algebra
over\/~$\bbQ(\sqrt{\Delta})(S)$.
It~corresponds to the restriction to the generic
point\/~$\smash{\eta \in S_{\bbQ(\sqrt{\Delta})}}$
of some lift\/
$\overline{c} \in \Br(S)$
of\/~$c$.\medskip

\noindent
{\bf Proof.}
{\em
We~have
$\smash{H^1(\Gal(\overline\bbQ/L), \Pic(S_{\overline\bbQ})) = 0}$
as the cubic surface
$S_L$
has Galois-invariant~sixers. The~inflation-restriction sequence yields
$$H^1(\Gal(\overline{\bbQ}/\bbQ(\sqrt{\Delta})), \Pic(S_{\overline\bbQ})) \cong H^1(\Gal(L/\bbQ(\sqrt{\Delta})), \Pic(S_L)) \, .$$
The~assertion easily~follows.
}
\eop
\end{lem}

\begin{remark}
It~is well known that a class
in~$\smash{\Br(S_{\bbQ(\sqrt{\Delta})})}$
is uniquely determined by its restriction
to~$\Br(\bbQ(\sqrt{\Delta})(S))$.
The~corresponding Azumaya algebra over the whole
of~$\smash{S_{\bbQ(\sqrt{\Delta})}}$
may be described as~follows.

Let~$x \in S_{\bbQ(\sqrt{\Delta})}$.
We~know that
$\div(\Psi)$
is the norm of a divisor
on~$S_L$.
That~one~is necessarily locally~principal. I.e.,~we have a rational function
$\smash{f_x}$
such that
$\div(N_{L/\bbQ(\sqrt{\Delta})} f_x) = \div(\Psi)$
on a Zariski neighbourhood
of~$x$.
Over~the maximal such
neighbourhood~$U_x$,
we define an Azumaya algebra~by
$$\smash{\textstyle Q_x := (\calO_{U_x} \!\otimes_{\bbQ(\sqrt{\Delta})}\! L)\{Y_x\} / (Y_x^3 - \frac\Psi{N_{L/\bbQ(\sqrt{\Delta})} f_x}) \, .}$$
Again,~we suppose
$Y_x t = \sigma(t) Y_x$
for~$t \in \calO_{U_x} \!\otimes_{\bbQ(\sqrt{\Delta})}\! L$.
In~particular, in a
neighbourhood~$U_\eta$
of the generic point, we have the Azumaya algebra
$Q_\eta := (\calO_{U_\eta} \!\otimes_{\bbQ(\sqrt{\Delta})}\! L) \{Y\} / (Y^3 - \Psi)$.

Over~$U_\eta \cap U_x$,
there is the isomorphism
$\iota_{\eta,x} \colon Q_\eta |_{U_\eta \cap U_x} \to Q_x |_{U_\eta \cap U_x}$,
given~by
$$Y \mapsto f_x Y_x \, .$$
For~two
points~$x, y \in S_{\bbQ(\sqrt{\Delta})}$,
the isomorphism
$\iota_{\eta,y} \!\circ\! \iota_{\eta,x}^{-1} \colon Q_x |_{U_\eta \cap U_x \cap U_y} \to Q_y |_{U_\eta \cap U_x \cap U_y}$
extends
to~$U_x \cap U_y$.
Hence,~the Azumaya algebras
$Q_x$
may be glued together along these~isomorphisms. This~yields an Azumaya algebra
$\calQ$
over~$S_{\bbQ(\sqrt{\Delta})}$.
\end{remark}

\begin{theo}
\label{Brexp}
Let\/~$\pi\colon S \to \Spec \bbQ$
be a non-singular cubic surface having a Galois invariant~triplet.
Further,~let\/
$c \in \Br(S)/\pi^*\!\Br(\bbQ)$
be a non-zero~class. Let,~finally,
$L$
be a splitting field that is cyclic of degree three over a
field\/~$\bbQ(\sqrt{\Delta})$
and\/
$\Psi$
a rational function
representing\/~$c$.\smallskip

\noindent
We~fix an isomorphism\/
$\Gal(L/\bbQ(\sqrt{\Delta})) \to \frac13\bbZ/\bbZ$.
For~every prime\/
$\frakp$
of\/~$\bbQ(\sqrt{\Delta})$
that does not split
in\/~$L$,
this fixes an isomorphism\/
$\smash{\Gal(L_\frakp/\bbQ(\sqrt{\Delta})_\frakp) \to \frac13\bbZ/\bbZ}$.\smallskip

\noindent
Then,~there is a
representative\/~$\underline{c} \in \Br(S)$
such that, for every prime
number\/~$p$,
the local evaluation~map\/
$\ev_p (\underline{c}, \,.\; ) \colon S(\bbQ_p) \to \bbQ/\bbZ$
is given as~follows.

\begin{abc}
\item
If~a prime\/
$\frakp$
above\/~$p$
splits in\/
$L/\bbQ(\sqrt{\Delta})$
then\/
$\ev_p (\underline{c}, x) = 0$
independently
of\/~$x$.
\item
Otherwise,~for\/
$x \in S(\bbQ_p) \subseteq S(\bbQ_p(\sqrt{\Delta}))$,
choose~a rational
function\/~$f_x$
on\/~$S_L$
such that, in a neighbourhood
of\/~$x$,
$\div(N\!f_x) = \div(\Psi)$.
Then,
$$\ev_p (\underline{c}, x) = 
\left\{
\begin{array}{cl}
\theta_\frakp(\Psi / N\!f_x) & {\it if~} (p) = \frakp {\it ~is~inert~in~} L, \\
2\theta_\frakp(\Psi / N\!f_x) & {\it if~} (p) = \frakp^2 {\it ~is~ramified~in~} L, \\
\theta_{\frakp_1}(\Psi / N\!f_x) + \theta_{\frakp_2}(\Psi / N\!f_x) & {\it if~} (p) = \frakp_1 \frakp_2 {\it ~splits~in~} L,
\end{array}
\right.$$
where\/~$\smash{\theta_\frakp \colon L_\frakp^* \to \Gal(L_\frakp/\bbQ(\sqrt{\Delta})_\frakp) \cong \frac13\bbZ/\bbZ}$
is the norm residue homomorphism for the extension\/
$L_\frakp/\bbQ(\sqrt{\Delta})_\frakp$
of local~fields.
\end{abc}

\noindent
{\bf Proof.}
{\em
The~assertion immediately follows from the~above.
}
\eop
\end{theo}

\begin{remso}[{\rm Diagonal cubic surfaces}{}]
\begin{iii}
\item
For~the particular case of a diagonal cubic surface
over~$\bbQ$,
the effect of the Brauer-Manin~obstruction has been studied intensively~\cite{CKS}. Our~approach covers diagonal cubic surfaces,~too.

Indeed,~on the surface given
by~$ax^3 + by^3 + cz^3 + dw^3 = 0$,
a Galois~invariant pair of Steiner trihedra is formed by the three planes given
by~$\sqrt[3]{a} x + \sqrt[3]{b} \zeta_3 y = 0$
together with the three planes given
by~$\sqrt[3]{c} z + \sqrt[3]{d} \zeta_3 w = 0$.
The~two other pairs completing the triplet are formed analogously pairing
$x$
with
$z$
or~$w$.
\item
Diagonal~cubic surfaces, are however, far from being the generic~case. There~are at least two~reasons. First,~the 27~lines are defined over
$\bbQ(\zeta_3, \sqrt[3]{a/d}, \sqrt[3]{b/d}, \sqrt[3]{c/d})$,
which is of degree at
most~$54$,
not~$216$.
Furthermore,~each of the six Steiner trihedra described has the property that the three planes intersect in a line while, generically, they intersect in a~point.
\item
For~the two pairs of Steiner trihedra defined by pairings
$(xy)(zw)$
and~$(xz)(yw)$,
we
find~$L = \bbQ(\zeta_3, \sqrt[3]{{ad}/{bc}})$
for the splitting field of the corresponding three double-sixes. This~is exactly the field occurring in~\cite[Lemme~1]{CKS}.

In~fact, the corresponding 18~lines break up into the six~trios
\begin{eqnarray*}
D_1^i\colon& \zeta_3^k \sqrt[3]{a/b} \,x + y = 0,\; \zeta_3^i z + \zeta_3^{-k} \sqrt[3]{d/c} \,w = 0, \quad k = 1,2,3 \,, \\ 
D_2^j\colon& \zeta_3^k \sqrt[3]{a/c} \,x + z = 0,\; \zeta_3^j y + \zeta_3^{-k}
\sqrt[3]{d/b} \,w = 0, \quad k = 1,2,3 \,.
\end{eqnarray*}
each consisting of pairwise skew~lines. Observe~that the trios are invariant under
$\Gal(\overline\bbQ/L)$
as
$\sqrt[3]{a/b}\sqrt[3]{d/c} = \sqrt[3]{a/c}\sqrt[3]{d/b} \in L$.

Further,~$D_1^i \cup D_2^j$
is a sixer if and only
if~$i \neq j$.
We~have the three double-sixes
$(D_1^1 \cup D_2^2) \cup (D_1^2 \cup D_2^1)$,
$(D_1^1 \cup D_2^3) \cup (D_1^3 \cup D_2^1)$,
and
$(D_1^2 \cup D_2^3) \cup (D_1^3 \cup D_2^2)$.
\item
Our~method to compute the Brauer-Manin~obstruction is, however, different from that of~\cite{CKS}. Colliot-Th\'el\`ene, Kanevsky, and Sansuc worked over larger extensions of the base~field.
\end{iii}
\end{remso}

\paragraph{\em Explicit Galois~descent.}

\begin{ttt}
Recall~that, in~\cite{EJ3}, we described a method to construct non-singular cubic surfaces
over~$\bbQ$
with a Galois invariant pair of Steiner~trihedra. The~idea was to start with cubic surfaces in Cayley-Salmon~form. For~these, we developed an explicit version of Galois~descent.
\end{ttt}

\begin{ttt}
More~concretely, given a quadratic number field
$D$
or~$D = \bbQ \oplus \bbQ$,
an element
$u \in D \!\setminus\! \{0\}$,
and a starting polynomial
$f \in D[T]$
of degree three without multiple zeroes, we constructed a cubic surface
$S^{u_0, u_1}_{(a_0, \ldots, a_5)}$
over~$\bbQ$
such~that
$$S^{u_0, u_1}_{(a_0, \ldots, a_5)} \times_{\Spec \bbQ} \Spec \overline\bbQ$$
is isomorphic to the surface
$S_{u_0, u_1}^{(a_0, \ldots, a_5)}$
in~$\bP^5_{\overline\bbQ}$
given~by
\begin{eqnarray*}
u_0 X_0 X_1 X_2 + u_1 X_3 X_4 X_5 & = & 0 \, , \\
a_0 X_0 + a_1 X_1 + a_2 X_2 + a_3 X_3 + a_4 X_4 + a_5 X_5 & = & 0 \, , \\
X_0 + \phantom{a_1} X_1 + \phantom{a_2} X_2 + \phantom{a_3} X_3 + \phantom{a_4} X_4 + \phantom{a_5} X_5 & = & 0 \, .
\end{eqnarray*}
Here,~$u_0 = \iota_0(u)$
and~$u_1 = \iota_1(u)$
for
$\iota_0, \iota_1 \colon D \to \overline\bbQ$
the two~embeddings
of~$D$.
Further,~$a \in A$
and
$a_i = \tau_i(a)$
for
$\tau_1, \ldots, \tau_5 \colon D[T]/(f) \to \overline\bbQ$
the six~embeddings of the \'etale
\mbox{$\bbQ$-algebra}
$D[T]/(f)$.
Thereby,~$\tau_0, \tau_1, \tau_2$
are supposed to be compatible
with~$\iota_0$,
whereas
$\tau_3, \tau_4, \tau_5$
are supposed to be compatible
with~$\iota_1$.
\end{ttt}

\begin{defi}
The~general cubic~polynomial
$$\Phi_{u_0, u_1}^{(a_0, \ldots, a_5)}(T) := \frac1{u_0} (a_0 + T)(a_1 + T)(a_2 + T) - \frac1{u_1} (a_3 + T)(a_4 + T)(a_5 + T)$$
is called the {\em auxiliary polynomial\/} associated
with~$S_{u_0, u_1}^{(a_0, \ldots, a_5)}$.
\end{defi}

\begin{prop}
\label{desc}
Let\/~$A$
be an \'etale
algebra of rank~three over a commutative semi\-simple\/
\mbox{$\bbQ$-}al\-ge\-bra\/~$D$
of dimension~two and\/
$u \in D$
as well
as~$a \in A$
as~above. Further,~suppose that\/
$S_{u_0, u_1}^{(a_0, \ldots, a_5)}$
is non-singular. Then,

\begin{abc}
\item
$\Phi_{u_0, u_1}^{(a_0, \ldots, a_5)} \in \bbQ[T]$
for\/
$D \cong \bbQ \oplus \bbQ$
and
$\sqrt{d}\,\Phi_{u_0, u_1}^{(a_0, \ldots, a_5)} \in \bbQ[T]$
for\/
$D \cong \bbQ(\sqrt{d})$.
\item
The~operation of an~element\/
$\sigma \in \Gal(\overline\bbQ/\bbQ)$
on\/~$S_{(a_0,\ldots,a_5)} \times_{\Spec \bbQ} \Spec \overline\bbQ$
goes over into the automorphism\/
$\smash{\pi_\sigma \!\circ\! t_\sigma \colon S^{(a_0, \ldots, a_5)} \to S^{(a_0, \ldots, a_5)} \, .}$

Here,~$\pi_\sigma$
permutes the coordinates according to the
rule\/~$a_{\pi_\sigma (i)} = \sigma(a_i)$
while\/
$t_\sigma$~is
the naive operation
of\/~$\sigma$
on\/~$S_{(a_0,\ldots,a_5)}$
as a morphism of schemes twisted
by\/~$\sigma$.
\item
Further,~on the descent
variety\/~$\smash{S^{u_0, u_1}_{(a_0, \ldots, a_5)}}$
over\/~$\bbQ$,
there~are
\begin{iii}
\item
nine~obvious lines given by
$L_{\{i,j\}} \colon \iota^* X_i = \iota^* X_j = 0$
for\/
$i = 0,1,2$
and\/~$j = 3,4,5$,
\item
18~non-obvious lines given~by\/
$\smash{L^\lambda_\rho \colon \iota^* Z_0 + \iota^* Z_{\rho(0)} = \iota^* Z_1 + \iota^* Z_{\rho(1)} = 0}$
for\/
$\lambda$
a zero of\/
$\smash{\Phi_{u_0, u_1}^{(a_0, \ldots, a_5)}}$
and\/~$\rho \colon \{0,1,2\} \to \{3,4,5\}$
a~bijection. Here,~the
coordinates\/~$Z_i$
are given~by
\begin{eqnarray*}
Z_0 := -Y_0 + Y_1 + Y_2, \quad
Z_1 := Y_0 - Y_1 + Y_2, \quad
Z_2 := Y_0 + Y_1 - Y_2, \\
Z_3 := -Y_3 + Y_4 + Y_5, \quad
Z_4 := Y_3 - Y_4 + Y_5, \quad
Z_5 := Y_3 + Y_4 - Y_5 \hspace{1.5mm}
\end{eqnarray*}
for\/~$Y_i := (a_i + \lambda) X_i, \quad i = 0,\ldots,5$.
\end{iii}

An~element\/
$\sigma \in \Gal(\overline\bbQ/\bbQ)$
acts on the lines according to the~rules
$$\sigma(L_{\{i,j\}}) = L_{\{\pi_\sigma(i), \pi_\sigma(j)\}}, \qquad \sigma(L^\lambda_\rho) = L^{\lambda^\sigma}_{\rho^{\pi_\sigma}} \, .$$
\item
The nine obvious lines are formed by a pair of Steiner~trihedra. The~twelve non-obvious lines corresponding to two zeroes of\/
$\smash{\Phi_{u_0, u_1}^{(a_0, \ldots, a_5)}}$
form a \mbox{double-six}.
\end{abc}\smallskip

\noindent
{\em
{\bf Proof.}
This~is a summary of results obtained in~\cite{EJ3}. a)~is \cite[Lemma~6.1]{EJ3}, while b) and~c) are \cite[Proposition~6.5]{EJ3}.
The~first assertion of~d) is \cite[Fact~4.2]{EJ3}. Finally,~the second one is shown \cite[Proposition~4.6]{EJ3} together with~\cite[Proposition~2.8]{EJ3}.
}
\eop
\end{prop}

\begin{rems}
\begin{iii}
\item
The~construction described is actually more general than needed~here. It~yields cubic surfaces with a Galois invariant pair of Steiner~trihedra. The~two other pairs, completing the first to a triplet, may be constructed as to be Galois invariant or to be defined over a quadratic number~field. Several~examples of both sorts were given in~\cite[Section~7]{EJ3}.
\item
Fix~the two Steiner trihedra corresponding to the 18 non-obvious~lines. Then,~parts c.ii) and~d) of Proposition~\ref{desc} together show that the field extension splitting the corresponding three double-sixes is the same as the splitting field of the auxiliary
polynomial~$\smash{\Phi_{u_0, u_1}^{(a_0, \ldots, a_5)}}$.\end{iii}
\end{rems}

\paragraph{\em Application: Manin's conjecture.}

\begin{ttt}
Recall~that a conjecture, due to Yu.~I.~Manin, asserts that the number of 
\mbox{$\bbQ$-ra}\-tio\-nal
points of anticanonical height
$\leq\! B$
on a Fano
variety~$S$
is asymptotically equal to
$\tau B \log^{\rk \Pic (S) - 1} B$,
for~$B\to\infty$.
Further,~the coefficient
$\tau \in \bbR$
is conjectured to be the Tamagawa-type
number~$\tau (S)$
introduced by E.~Peyre in~\cite{Pe}. In~the particular case of a cubic surface, the anticanonical height is the same as the naive~height.
\end{ttt}

\begin{ttt}
E.~Peyre's Tamagawa-type~number is defined in \cite[Definition~2.4]{PT}~as
$$\smash{\tau(S) := \alpha(S) \!\cdot\! \beta(S) \cdot \lim_{s\to1} \, (s-1)^t L(s,\chi_{\Pic(S_{\overline\bbQ})}) \cdot \tau_H \!\big( S(\bbA_\bbQ)^{\rm Br} \big)}$$
for~$t = \rk \Pic(S)$.
%
%
%
Here,~$\tau_H$
is the Tamagawa measure on the set
$S(\bbA_\bbQ)$
of adelic points
on~$S$
and
$S(\bbA_\bbQ)^{\rm Br} \subseteq S(\bbA_\bbQ)$~denotes
the part which is not affected by the Brauer-Manin~obstruction. For~details, in particular on the constant factors, we refer to the original~literature. But~observe that the Tamagawa measure is the product measure of measures
on~$S(\bbQ_\nu)$
for~$\nu$
running through the places
of~$\bbQ$.
\end{ttt}

%
%

\begin{ttt}
Using~\cite[Algorithm~5.8]{EJ3}, we constructed many examples of smooth cubic surfaces
over~$\bbQ$
with a Galois~invariant pair of Steiner~trihedra. For~each of them, one may apply Theorem~\ref{Brexp} to compute the effect of the Brauer-Manin~obstruction. Then,~the method described in~\cite{EJ1} applies for the computation of Peyre's~constant.

From~the considerable supply, the examples below were chosen in the hope that they indicate the main~phenomena. The~Brauer-Manin~obstruction may work at many primes simultaneously but examples where few primes are involved are more~interesting. We~will show that the fraction of the Tamagawa measure excluded by the obstruction can vary~greatly.
\end{ttt}

\begin{ex}
Start with
$g(U) := U^2 - 1$,
i.e.~$D = \bbQ \oplus \bbQ$,
$u_0 := 1$,
$u_1 := \frac12$,
$f_0(V) := V^3 + 6V^2 + 9V + 1$,
and~$f_1(V) := V^3 + \frac92V^2 + \frac92V$.
Then,~the auxiliary polynomial
is~$-V^3 - 3V^2 + 1$.
\cite[Algorithm~5.8]{EJ3} yields the cubic
surface~$S$
given by the equation
$$T_0^3 - T_0^2T_2 - T_0^2T_3 - 2T_0T_2^2 + T_0T_2T_3 - T_1^3 + 3T_1^2T_2 - 3T_1T_2T_3 + 3T_1T_3^2 - T_2^3 - T_2^2T_3 + T_3^3 = 0 \, .$$
This~example is constructed in such a way that
$f_1$
completely splits, while
$f_0$
and the auxiliary polynomial have the same splitting field that is an
$A_3$-extension
of~$\bbQ$.
Hence,~the Galois group operating on the 27~lines is of order~three. We~have orbit structure
$[3,3,3,3,3,3,3,3,3]$
and~$\smash{H^1(\Gal(\overline{\bbQ}/\bbQ), \Pic(S_{\overline\bbQ})) \cong \bbZ/3\bbZ \times \bbZ/3\bbZ}$.

$S$~has
bad reduction at the primes
$3$
and~$19$.
It~turns out that the local evaluation map
at~$19$
is~constant. Hence,~the Brauer-Manin obstruction works only at the
prime~$3$.

Here,~the reduction is the cone over a cubic curve
over~$\bbF_{\!3}$
having a~cusp. Only~the nine smooth points lift to
\mbox{$3$-adic}~ones.
Among~them, exactly one is allowed by the Brauer-Manin~obstruction. Hence,~from the whole
of~$S(\bbQ_3)$,
which is of
measure~$1$,
only a subset of
measure~$\frac19$
is~allowed.

Using~this, for Peyre's constant, we find
$\tau(S) \approx 0.1567$.
There~are actually 599
$\bbQ$-rational
points of height at most
$4000$
in comparison with a prediction
of~$627$.
\end{ex}

\begin{remark}
It~is, may be, of interest to note that
$S$
has good reduction
at~$2$
and the reduction has exactly one
\mbox{$\bbF_{\!2}$-rational}~point.
An~example of such a cubic surface
over~$\bbF_{\!2}$
was once much sought-after and finally found by Swinnerton-Dyer~\cite{SD2}.
\end{remark}

\begin{ex}
\label{bspzwei}
Start with
$g(U) := U^2 - 1$,
i.e.~$D = \bbQ \oplus \bbQ$,
$u_0 := -\frac13$,
$u_1 := 1$,
$f_0(V) := V^3 - V + \frac13$,
and~$f_1(V) := V^3 + 3V^2 - 4V + 1$.
Then,~the auxiliary polynomial
is~$-4V^3 - 3V^2 + 7V - 2$.
\cite[Algorithm~5.8]{EJ3} yields the cubic
surface~$S$
given by the~equation
\begin{eqnarray*}
&& -3T_0^3 - 6T_0^2T_1 - 3T_0^2T_2 + 3T_0^2T_3 - 3T_0T_1^2 + 3T_0T_1T_3 + 3T_0T_2^2 + 6T_0T_3^2 + 2T_1^3 \\
&& {}\!\! - 4T_1^2T_2 - T_1^2T_3 + 10T_1T_2^2 - 4T_1T_2T_3 - 9T_1T_3^2 + 6T_2^3 - 8T_2^2T_3 - 8T_2T_3^2 + 4T_3^3 = 0 \, .
\end{eqnarray*}
Here,
$f_0$,
$f_1$,
and the auxiliary polynomial have splitting fields that are linearly disjoint
$A_3$-extensions
of~$\bbQ$.
Hence,~the Galois group operating on the 27~lines is of
order~$27$.
It~is the maximal
$3$-group
stabilizing a pair of Steiner~trihedra. We~have orbit structure
$[9,9,9]$
and~$\smash{H^1(\Gal(\overline{\bbQ}/\bbQ), \Pic(S_{\overline\bbQ})) \cong \bbZ/3\bbZ}$.

The~three minimal splitting fields are given by the polynomials
$T^3 + T^2 - 10T - 8$,
$T^3 - 21T + 28$,
and~$T^3 - 21T - 35$.
Observe~that the first one is unramified at both primes,
$3$
and~$7$,
while the other two fields~ramify.

$S$~has
bad reduction at the primes
$3$,
$7$,
and~$31$.
The~local
\mbox{$H^1$-criterion}~\cite[Lemma~6.6]{EJ2}
excludes the
prime~$31$.
Hence,~the Brauer-Manin obstruction works only at the primes
$3$
and~$7$.

At~$3$,
the reduction is the cone over a cubic curve
over~$\bbF_{\!3}$
having a~cusp. The~only singular point that lifts to a
\mbox{$3$-adic}~one
is~$(1:0:1:1)$.
The~local evaluation map is constant on the points of smooth reduction and decomposes the others into two equal~parts. This~means,
$S(\bbQ_3)$
breaks into sets of measures
$1$,
$\frac13$,
and
$\frac13$,~respectively.

At~$7$,
the reduction is Cayley's ruled cubic surface
over~$\bbF_{\!7}$.
Exactly~one singular point lifts to a
\mbox{$7$-adic}~one.
This~is~$(1:5:1:0)$.
Again,~the local evaluation map is constant on the points of smooth reduction and decomposes the others into two equal~parts.
Thus,~$S(\bbQ_7)$
breaks into sets of measures
$1$,
$\frac17$,
and
$\frac17$.
An~easy calculation shows
$\tau_H(S(\bbA_\bbQ)^{\Br}) = \frac{11}{45} \tau_H(S(\bbA_\bbQ))$.

Using~this, for Peyre's constant, we find
$\tau(S) \approx 1.7311$.
There~are actually 6880
$\bbQ$-rational
points of height at most
$4000$
in comparison with a prediction
of~$6924$.
\end{ex}

\begin{ex}
Start with
$g(U) := U^2 - 1$,
i.e.~$D = \bbQ \oplus \bbQ$,
$u_0 := \frac{13}7$,
$u_1 := -\frac{31}{169}$,
$f_0(V) := V^3 - 16V^2 + 85V - \frac{1049}7$,
and~$f_1(V) := V^3 - \frac{337}{31}V^2 + \frac{1216}{31}V - 47$.
Then,~the auxiliary polynomial
is~$\frac{2414}{403}V^3 - \frac{848021}{12493}V^2 + \frac{595}{13}V - \frac{1537566}{12493}$.
\cite[Algorithm~5.8]{EJ3} yields the cubic
surface~$S$
given by the~equation
\begin{eqnarray*}
&& 13T_0^3 - 8T_0^2T_1 + 9T_0^2T_2 + 44T_0^2T_3 - 9T_0T_1^2 - 5T_0T_1T_2 + T_0T_1T_3 + 4T_0T_2^2 \\
&& \hspace{0.2cm} {} - 19T_0T_2T_3 - 61T_0T_3^2 - T_1^2T_2 - 24T_1^2T_3 + 3T_1T_2^2 + 42T_1T_2T_3 - 16T_1T_3^2 \\
&& \hspace{6.2cm} {} + 2T_2^3 + 10T_2^2T_3 - 60T_2T_3^2 + 6T_3^3 = 0 \, .
\end{eqnarray*}
Here,~$f_0$,
and~$f_1$
define
$A_3$-extensions
of~$\bbQ$,
ramified only at
$7$
and~$13$,~respectively.
As~the auxiliary polynomial is of
type~$S_3$,
the Galois group operating on the 27~lines is of
order~$54$.
We~have orbit structure
$[9,9,9]$
and~$\smash{H^1(\Gal(\overline{\bbQ}/\bbQ), \Pic(S_{\overline\bbQ})) \cong \bbZ/3\bbZ}$.

The~minimal splitting fields are given by the polynomials
$T^3 - T^2 - 30T - 27$,
$T^3 - T^2 - 30T + 64$,
and~$T^3 - 75\,986\,365T - 753\,805\,852\,436$.
The~first two are cyclic of degree~three. Both~are ramified exactly at
$7$
and~$13$.
The~third splitting field is of
type~$S_3$
and ramified at
$7$,
$13$,
$127$,
and~$387\,512\,500\,241$.

$S$~has
bad reduction at the primes
$3$,
$7$,
$13$,
$127$,
$281$,
$84\,629$
and~$387\,512\,500\,241$.
The~local
\mbox{$H^1$-criterion}
excludes all primes except
$7$
and~$13$.
Hence,~the Brauer-Manin obstruction works only at the primes
$7$
and~$13$.

At~$7$,
the reduction is a normal cubic surface with three singularities of
type~$A_2$
that are defined
over~$\bbF_{\!7^3}$.
The~local evaluation map
on~$S(\bbQ_7)$
factors
via~$S(\bbF_{\!7})$
and decomposes these 57~points into three equal~classes. This~means,
$S(\bbQ_7)$
breaks into three sets, each of which is of
measure~$\frac{19}{49}$.

At~$13$,
the reduction is a normal cubic surface having three singularities of
type~$A_2$
that are defined
over~$\bbF_{\!13}$.
The~behaviour is a bit more complicated as there are
\mbox{$13$-adic}~points
reducing to the~singularities. On~the points having good reduction, the local evaluation map factors
via~$S(\bbF_{\!13})$.
It~splits the 180 smooth points into three classes of 60~elements~each.

Summarizing,~we clearly have
$\tau_H(S(\bbA_\bbQ)^{\Br}) = \frac13 \tau_H(S(\bbA_\bbQ))$.
Using~this, for Peyre's constant, we find
$\tau(S) \approx 0.5907$.
There~are actually 2370
$\bbQ$-rational
points of height at most
$4000$
in comparison with a prediction
of~$2363$.
\end{ex}

\begin{remark}
On~$S(\bbF_{\!7})$
and on the smooth
$\bbF_{\!13}$-rational
points, we find exactly the decompositions defined by Mordell-Weil~equivalence~\cite[Corollary~3.4.3]{EJ4}.
\end{remark}

\begin{ex}
\label{bspvier}
Start with
$g(U) := U^2 - 7$,
i.e.~$D = \bbQ(\sqrt{7})$,
$\smash{u := \frac1{-1 + 2\sqrt{7}} = -\frac{1 + 2\sqrt{7}}{29}}$,
and
$$f(V) := (-1 + 2\sqrt{7}) V^3 + (1 + 2\sqrt{7}) V^2 + (-2 + \sqrt{7}) V + (-1 + 2\sqrt{7}) \, .$$
Then,~the auxiliary polynomial is
$4\sqrt{7} (V^3 + V^2 + \frac12V + 1)$.
\cite[Algorithm~5.8]{EJ3} yields the cubic
surface~$S$
given by the~equation
\begin{eqnarray*}
&& 11T_0^3 + 32T_0^2T_1 + 31T_0^2T_2 + 52T_0^2T_3 - 33T_0T_1^2 - 93T_0T_1T_2 - 36T_0T_1T_3 + 37T_0T_2^2 \\
&& {} + 46T_0T_2T_3 - 34T_0T_3^2 + 22T_1^3 -  10T_1^2T_2 - 21T_1^2T_3 + 75T_1T_2^2 - 4T_1T_2T_3 + 30T_1T_3^2 \\
&& \hspace{7.5cm} {} + 133T_2^3 + 34T_2^2T_3 - 8T_2T_3^2 + 2T_3^3 = 0 \, .
\end{eqnarray*}
Here,
$f$
has Galois
group~$S_3$.
Its~discriminant
is~$(-22144 + 3806\sqrt{7})$,
a number of
norm~$19722^2$.
This~ensures that not only one, but three pairs of Steiner trihedra are Galois~invariant and that we have orbit
structure~$[9,9,9]$.
The~auxiliary polynomial is of
type~$S_3$.
Hence,~the Galois group operating on the 27~lines is the maximal possible
group~$U_t$
stabilizing a~triplet. It~is of
order~$216$.

The~three minimal splitting fields are given by the polynomials
$T^3 - 7T - 8$,
$T^3 + T^2 - 177T - 2059$,
and~$T^3 - T^2 + 32T + 368$.

$S$~has
bad reduction at the primes
$2$,
$7$,
$19$,
$89$,
$151$
and~$173$.
The~local
\mbox{$H^1$-criterion}
excludes
$2$,
$7$,
$89$,
$151$
and~$173$.
Hence,~the Brauer-Manin obstruction works only at the
prime~$19$.

Here,~the reduction is to three planes
over~$\bbF_{\!19}$.
The~intersection points do not lift to
\mbox{$\bbQ_{19}$-rational}~ones.
The~Brauer-Manin obstruction allows exactly one of the three~planes. This~means, from the whole
of~$S(\bbQ_{19})$,
which is of
measure~$3$,
only a subset of
measure~$1$
is~allowed.

Using~this, for Peyre's constant, we find
$\tau(S) \approx 0.0553$.
There~are actually 216
$\bbQ$-rational
points of height at most
$4000$
in comparison with a prediction
of~$221$.
\end{ex}

\begin{rems}
\begin{iii}
\item
Except for the first, these examples have orbit
structures~$[9,9,9]$.
Unfortunately,~the rational functions representing the Brauer class are not very well suited for a~reproduction. For~instance, in Example~\ref{bspzwei}, the non-trivial Brauer class can be represented by the cubic~form
\begin{eqnarray*}
21965T_0^3 + 55863T_0^2T_2 - 53607T_0^2T_3 - 1215T_0T_1^2 + 75402T_0T_1T_2 - 136125T_0T_1T_3 && \\
{} + 35961T_0T_2^2 - 133200T_0T_2T_3 + 44382T_0T_3^2 + 8402T_1^3 + 149304T_1^2T_2 && \\
{} - 272189T_1^2T_3 - 22210T_1T_2^2 - 249626T_1T_2T_3 + 526313T_1T_3^2 - 36518T_2^3 && \\
{} - 70576T_2^2T_3 + 184098T_2T_3^2 - 26618T_3^3 &&
\end{eqnarray*}
divided by the cube of a linear~form.
\item
In~the final example, the rational functions representing the Brauer class have coefficients
in~$\bbQ(\sqrt{7})$.
Rational~functions with coefficients
in~$\bbQ$
do not serve our~purposes.

Indeed,~the prime
$19 = \frakp_1 \frakp_2$
splits
in~$\bbQ(\sqrt{7})$.
According~to class field theory, the two norm residue homomorphisms
$\theta_{\frakp_1}, \theta_{\frakp_2}$
are opposite to each other when applied to elements
from~$\bbQ_{19}$.
Thus,~$\Psi \in \bbQ[T_0, \ldots, T_3]$
would imply that the Brauer-Manin~obstruction does not work
at~$19$,
at least not on
\mbox{$\bbQ$-rational}~points.
\item
To~apply the local
\mbox{$H^1$-criterion},
we compute a
degree-$36$
resolvent the zeroes of which correspond to the double-sixes. According~to Lemma~\ref{doppelsechs}, we have
$\smash{H^1 (\Gal(\overline\bbQ_p/\bbQ_p), \Pic(S_{\overline\bbQ})) = 0}$
when the resolvent has a
\mbox{$p$-adic}~zero.
\item
It~is noticeable from the examples that the reduction types at the relevant primes are distributed in an unusual~way. This~is partially explained by \cite[Proposition~1.12]{EJ4}.
\end{iii}
\end{rems}

\section{Two triplets}

\begin{lem}
Let\/~$\calT = \{T_1, T_2, T_3\}$
and\/~$\calT^\prime = \{T_1^\prime, T_2^\prime, T_3^\prime\}$
be two different decompositions of the 27~lines into pairs of Steiner~trihedra. Then,~up to permutation of rows and columns, for the
matrix\/~$J$
defined
by\/~$J_{ij} := \#\!I_{ij}$
for\/~$I_{ij} := T_i \cap T_j^\prime$,
there are two~possibilities.

\begin{ABC}
\item
$
\left(
\begin{array}{ccc}
3 & 3 & 3 \\
3 & 3 & 3 \\
3 & 3 & 3
\end{array}
\right),
$
\item
$
\left(
\begin{array}{ccc}
5 & 2 & 2 \\
2 & 5 & 2 \\
2 & 2 & 5
\end{array}
\right).
$
\end{ABC}\smallskip

\noindent
{\bf Proof.}
{\em
As~$W(E_6)$
transitively operates on decompositions, we may suppose that
$\calT = \St_{(123)(456)}$
is the standard~decomposition.

Then,~the other nine decompositions of
type~$\St_{(ijk)(lmn)}$
yield type~B. Among~the 30~decompositions of
type~$\St_{(ij)(kl)(mn)}$,
there are exactly~18 for which one of the sets
$\{i,j\}$,
$\{k,l\}$,
or~$\{m,n\}$
is a subset
of~$\{1,2,3\}$.
These~lead to type~B,~too. The~twelve remaining triplets yield type~A.
}
\eop
\end{lem}

\begin{rems}
\begin{iii}
\item
Correspondingly, the group
$U_t \subset W(E_6)$
that stabilizes a triplet, operates on the set of all decompositions such that the orbits have lengths
$1$,
$12$,
and~$27$.
\item
The~group
$U_t$
operates as well on the set of the 240~triplets. Here,~the orbit lengths are
$1$, $1$, $1$, $1$, $1$, $1$,
$27$, $27$, $27$, $27$, $27$, $27$,
and~$72$.
In~fact, in case~A,
$U_t$
is able to permute the pairs
$T_1^\prime$,
$T_2^\prime$,
and~$T_3^\prime$
although
$T_1$,
$T_2$,
and~$T_3$
remain~fixed. The~same cannot happen in case~B.
\end{iii}
\end{rems}

\begin{exo}[{\rm Type~A}{}]
A~representative case for this type is given by the decomposition
$\calT = \St_{(123)(456)}$,
which, by definition, consists of the pairs
$$\left[
\begin{array}{ccc}
a_1 & b_2 & c_{12} \\
b_3 & c_{23} & a_2 \\
c_{13} & a_3 & b_1
\end{array}
\right],
\qquad
\left[
\begin{array}{ccc}
a_4 & b_5 & c_{45} \\
b_6 & c_{56} & a_5 \\
c_{46} & a_6 & b_4
\end{array}
\right],
\quad
{\rm and}
\quad
\left[
\begin{array}{ccc}
c_{14} & c_{15} & c_{16} \\
c_{24} & c_{25} & c_{26} \\
c_{34} & c_{35} & c_{36}
\end{array}
\right]
$$
of Steiner trihedra and
$\calT^\prime = \St_{(14)(25)(36)}$,
consisting of
$$\left[
\begin{array}{ccc}
a_1 & b_5 & c_{15} \\
b_2 & a_4 & c_{24} \\
c_{12} & c_{45} & c_{36}
\end{array}
\right],
\qquad
\left[
\begin{array}{ccc}
a_2 & b_6 & c_{26} \\
b_3 & a_5 & c_{35} \\
c_{23} & c_{56} & c_{14}
\end{array}
\right],
\quad
{\rm and}
\quad
\left[
\begin{array}{ccc}
a_3 & b_4 & c_{34} \\
b_1 & a_6 & c_{16} \\
c_{13} & c_{46} & c_{25}
\end{array}
\right].
$$
For~the matrix of the mutual intersections, one~obtains
$$
I =
\left(
\begin{array}{ccc}
\{a_1,b_2,c_{12}\} & \{a_2,b_3,c_{23}\} & \{a_3,b_1,c_{13}\} \\
\{a_4,b_5,c_{45}\} & \{a_5,b_6,c_{56}\} & \{a_6,b_4,c_{46}\} \\
\{c_{15},c_{24},c_{36}\} & \{c_{14},c_{26},c_{35}\} & \{c_{16},c_{25},c_{34}\}
\end{array}
\right).
$$
The~maximal subgroup
$U_{tt}$
of~$W(E_6)$
stabilizing both triplets is of
order~$216/72 = 3$.
It~is clearly the same as the maximal subgroup stabilizing each of the nine sets of size three obtained as~intersections.
\end{exo}

\begin{remark}
$U_{tt}$~actually
stabilizes four~decompositions of the 27~lines into three sets of~nine. They~are given by the rows, the columns as well as the positive and negative diagonals of the Sarrus~scheme. 
\end{remark}

\begin{remark}
Each~decomposition gives rise to eight {\em enneahedra,} systems of nine planes containing all the 27~lines. Two~decompositions meeting of type~A have an enneahedron in~common, as one may see at the
matrix~$I$.
Hence,~there are
$\frac{40 \cdot 12}2 / (\atop{4}{2}) = 40$
enneahedra, each of which originates from four~decompositions. On~the other hand, the remaining 160~enneahedra each stem from a single~one.

In~this form, the two types, in which two decompositions may be correlated, were known to L.\,Cremona~\cite{Cr} in~1870. Cremona~calls the~40 enneahedra of the first kind and the others of the~second.
\end{remark}

\begin{theo}
Let\/~$\pi\colon S \to \Spec\bbQ$
be a non-singular cubic~surface. Suppose that\/
$S$
has two\/
$\Gal(\overline\bbQ/\bbQ)$-invariant
triplets that do not correspond to the same decomposition and an orbit structure of\/
$[3,3,3,3,3,3,3,3,3]$
on the 27~lines.

\begin{abc}
\item
Then,~$\Br(S)/\pi^*\!\Br(\bbQ) \cong \bbZ/3\bbZ \times \bbZ/3\bbZ$.
\item
There~are actually four\/
$\Gal(\overline\bbQ/\bbQ)$-invariant
decompositions
on\/~$S$.
The~eight non-zero elements
of\/~$\Br(S)/\pi^*\!\Br(\bbQ)$
are the images of these under the class~map. In~other words,
$\cl$
induces a~bijection
$$\Theta_S^{\Gal(\overline\bbQ/\bbQ)} / A_3 \stackrel{\cong}{\longrightarrow} (\Br(S)/\pi^*\!\Br(\bbQ)) \setminus \{0\} \, .$$
\end{abc}
{\bf Proof.}
{\em
Observe~that,
on~$S$,
there is an enneahedron consisting of Galois invariant~planes. This~immediately implies the first assertion of~b).\medskip

\noindent
a)
Once~again, we will use the isomorphism
$$\delta\colon H^1(\Gal(\overline{\bbQ}/\bbQ), \Pic(S_{\overline\bbQ})) \to \Br(S) / \pi^*\!\Br(\bbQ) \, .$$
Further,~by Shapiro's lemma,
$\smash{H^1(\Gal(\overline{\bbQ}/\bbQ), \Pic(S_{\overline\bbQ})) \cong H^1(G, \Pic(S_{\overline\bbQ}))}$
for~$G$
the subgroup
of~$W(E_6)$,
through which
$\Gal(\overline{\bbQ}/\bbQ)$
actually operates on the 27~lines.

The~assumption clearly implies that the two decompositions considered meet each other of type~A.
Hence,~up to conjugation,
$G \subseteq U_{tt}$.
As~$U_{tt}$
is of order three, there are no proper subgroups, except for the trivial~group. Hence,~$G = U_{tt}$.
It~is now a direct calculation to verify
$\smash{H^1(U_{tt}, \Pic(S_{\overline\bbQ})) = \bbZ/3\bbZ \times \bbZ/3\bbZ}$.
This~was worked out by Yu.\,I.~Manin in~\cite[Lemma~45.4.1]{Ma}.\medskip\smallskip

\noindent
b)
We~will show this assertion in several~steps.\medskip

\noindent
{\em First step.}
Divisors.\smallskip

\noindent
We~will use Manin's~formula. The~orbits of size three define nine divisors, which we call
$\Delta_{ij}$
for~$1 \leq i,j \leq 3$.
Thereby,~$\Delta_{ij}$
shall correspond to the set
$I_{ij}$
for~$I$
the matrix of intersections given~above.

The~norm of a line is always one of
the~$\Delta_{ij}$.
Hence,~$N\!D = \{\sum_{ij} n_{ij} \Delta_{ij} \mid n_{ij} \in \bbZ \}$.
Further,~a direct calculation shows that
$L\Delta_{ij} = 1$
for every
line~$L$
on~$\smash{\Pic(S_{\overline\bbQ})}$
and every
pair~$(i,j)$.
Therefore,~a
divisor~$\sum_{ij} n_{ij} \Delta_{ij}$
is principal if and only
if~$\sum_{ij} n_{ij} = 0$.

Now,~consider a divisor of the form
$$\Delta_{i_1j_1} + \Delta_{i_2j_2} + \Delta_{i_3j_3}$$
such that
$$i_1 + i_2 + i_3 \equiv j_1 + j_2 + j_3 \equiv 0 \pmod 3 \, .$$
It~is not hard to check that such a divisor is always the norm of the sum of three lines in a~plane. Consequently,~the difference of two such divisors is
in~$N\!D_0$.\medskip\pagebreak[3]

\noindent
{\em Second step.}
The~nine residue~classes.\smallskip

\noindent
We~therefore have the following parallelogram~rule.
$$\Delta_{ij} + \Delta_{i^\prime j^\prime} - \Delta_{i+i^\prime, j+j^\prime} - \Delta_{33} = 0 \in (N\!D \cap D_0) / N\!D_0 \, ,$$
where the arithmetic of the indices is done
modulo~$3$.
Indeed,~$\Delta_{ij} + \Delta_{i^\prime j^\prime} + \Delta_{-i-i^\prime, -j-j^\prime}$
and~$\Delta_{i+i^\prime, j+j^\prime} + \Delta_{33} + \Delta_{-i-i^\prime, -j-j^\prime}$
are norms of sums of three lines in a~plane.

The~parallelogram rule implies, in particular, that every element of
$(N\!D \cap D_0) / N\!D_0$
is equal to some
$\Delta_{ij} - \Delta_{33}$.
Clearly,~these nine elements are mutually~distinct.
Further,~$(\Delta_{ij} - \Delta_{33}) + (\Delta_{i^\prime j^\prime} - \Delta_{33}) = \Delta_{i+i^\prime, j+j^\prime} - \Delta_{33}$.\medskip

\noindent
{\em Third step.}
A~non-zero Brauer~class.\smallskip

\noindent
Let~now a non-zero homomorphism
$g\colon (N\!D \cap D_0) / N\!D_0 \to \frac13\bbZ/\bbZ$
be~given. Then,
$\ker g = \{0, \Delta_{ij} - \Delta_{33} , \Delta_{2i,2j} - \Delta_{33}\}$
for a suitable
pair~$(i,j)$.
This~shows~that
$$\textstyle g^{-1}(\frac13) = \{\Delta_{i_1,j_1} - \Delta_{33}, \Delta_{i_2,j_2} - \Delta_{33}, \Delta_{i_3,j_3} - \Delta_{33}\}$$
for~$i_1 + i_2 + i_3 \equiv j_1 + j_2 + j_3 \equiv 0 \pmod 3$.

The~divisor
$\Delta_{i_1,j_1} + \Delta_{i_2,j_2} + \Delta_{i_3,j_3}$
consists exactly of the nine lines defined by a
pair~$T_1$
of Steiner~trihedra that is fixed
by~$U_{tt}$.
Further,
$\Delta_{33} + \Delta_{ij} + \Delta_{2i,2j}$
consists of the lines defined by another
pair~$T_2$.
Let~$\calT = (T_1,T_2,T_3)$
be the corresponding~triplet. We~claim that the restriction
of~$\cl(\calT)$
is exactly the cohomology class represented
by~$g$.

For~this, observe that the restriction map in group cohomology is dual to the norm map
$N_G\colon (N\!D \cap D_0) / N\!D_0 \to (N_G D \cap D_0) / N_G D_0$
for~$G$
any group causing the orbit
structure~$[9,9,9]$
given by the
triplet~$\calT$.
This~norm map sends
$\Delta_{i_1,j_1}$,
$\Delta_{i_2,j_2}$,
and~$\Delta_{i_3,j_3}$
to~$D_2$
and~$\Delta_{33}$
to~$D_1$.
I.e.,~$(\Delta_{i_l,j_l} - \Delta_{33})$
is mapped
to~$(D_2 - D_1)$
in the notation of~\ref{cong}.
As~$\cl(\calT)$
sends
$(D_2 - D_1)$
to~$\frac13$,
this is enough to imply the~claim.\medskip

\noindent
{\em Fourth step.}
Conclusion.\smallskip

\noindent
Consequently,~the mapping
$\smash{\Psi_S^{\Gal(\overline\bbQ/\bbQ)} / A_3 \to (\Br(S)/\pi^*\!\Br(\bbQ)) \setminus \{0\}}$
induced
by~$\cl$
is a~surjection. As,~on both sides, the sets are of size eight, bijectivity~follows.
}
\eop
\end{theo}

\begin{remark}
Up~to conjugation,
$G = U_{tt}$
is the only subgroup
of~$W(E_6)$
such that
$H^1(G, \Pic(S_{\overline\bbQ})) \cong \bbZ/3\bbZ \times \bbZ/3\bbZ$.
An~example of a cubic surface
over~$\bbQ$,
for which
$U_{tt}$
appears as the Galois group operating on the 27~lines, was given in~\cite[Example~7.6]{EJ3}.
Over~$\bbQ(\zeta_3)$,
one may simply consider the surface given
by~$x^3 + y^3 + z^3 + aw^3 = 0$
for~$a \in \bbQ^*$
a \mbox{non-cube}.
\end{remark}

\begin{exo}[{\rm Type~B}{}]
A~representative case for this type is given by the decompositions
$\calT = \St_{(123)(456)}$,
consisting of the pairs
$$\left[
\begin{array}{ccc}
a_1 & b_2 & c_{12} \\
b_3 & c_{23} & a_2 \\
c_{13} & a_3 & b_1
\end{array}
\right],
\qquad
\left[
\begin{array}{ccc}
a_4 & b_5 & c_{45} \\
b_6 & c_{56} & a_5 \\
c_{46} & a_6 & b_4
\end{array}
\right],
\qquad
\left[
\begin{array}{ccc}
c_{14} & c_{15} & c_{16} \\
c_{24} & c_{25} & c_{26} \\
c_{34} & c_{35} & c_{36}
\end{array}
\right]
$$
of Steiner trihedra, and
$\calT^\prime = \St_{(12)(34)(56)}$
consisting of
$$\left[
\begin{array}{ccc}
a_1 & b_4 & c_{14} \\
b_3 & a_2 & c_{23} \\
c_{13} & c_{24} & c_{56}
\end{array}
\right],
\qquad
\left[
\begin{array}{ccc}
a_3 & b_6 & c_{36} \\
b_5 & a_4 & c_{45} \\
c_{35} & c_{46} & c_{12}
\end{array}
\right],
\qquad
\left[
\begin{array}{ccc}
a_5 & b_2 & c_{25} \\
b_1 & a_6 & c_{16} \\
c_{15} & c_{26} & c_{34}
\end{array}
\right].
$$
For~the matrix of the mutual intersections, we~obtain
$$
\left(
\begin{array}{ccc}
\{a_1,a_2,b_3,c_{13},c_{23}\} & \{a_3,c_{12}\} & \{b_1,b_2\} \\
\{b_4,c_{56}\} & \{a_4,b_5,b_6,c_{45},c_{46}\} & \{a_5,a_6\} \\
\{c_{14},c_{24}\} & \{c_{35},c_{36}\} & \{c_{15},c_{16},c_{25},c_{26},c_{34}\}
\end{array}
\right).
$$
The~maximal subgroup
$U_{tt}^\prime \subset W(E_6)$
stabilizing both triplets is of
order~$216/27 = 8$.
The~Brauer classes defined by both triplets are annihilated when restricted to a surface on which these six pairs of trihedra are Galois~invariant.
\end{exo}

\begin{remark}
Observe~that
$U_{tt}^\prime$
does not stabilize any pair of Steiner trihedra different from the six given~above.
\end{remark}

\end{document}